\documentclass[12pt]{amsart}
\usepackage[margin=2.95cm]{geometry}

\usepackage{amssymb}
\usepackage{eucal}

\usepackage[all,cmtip,arrow,2cell,poly,curve]{xy}
\UseTwocells
\everyxy={0;<8pt,0pt>:;}

\newcommand{\bbZ}{{\mathbb{Z}}}
\newcommand{\bbC}{{\mathbb{C}}}

\newcommand{\bbN}{{\mathbb{N}}}
\newcommand{\gT}{\mathbb{T}^\lambda}
\newcommand{\gX}{\mathbb{X}^\lambda}

\newcommand{\calC}{{\mathcal{C}}}

\newcommand{\calX}{{\mathcal{X}}}

\newcommand{\op}{{\mathrm{op}}}
\newcommand{\catX}{{\mathcal{X}}}
\newcommand{\catY}{{\mathcal{Y}}}
\newcommand{\catZ}{{\mathcal{Z}}}
\newcommand{\catB}{{\mathcal{B}}}
\newcommand{\catC}{{\mathcal{C}}}
\newcommand{\catD}{{\mathcal{D}}}

\newcommand{\Stack}[1]{{\underline{#1}}} 

\newcommand{\GrpdBibu}{{\mathcal{G}\mathrm{rpd}\mathcal{B}\mathrm{ibu}}}
\newcommand{\LieGrpdPrBibu}{{\mathrm{st}\mathcal{M}\mathrm{fld}}}
\newcommand{\stMfld}{{\LieGrpdPrBibu}}
\newcommand{\Mfld}{{\mathcal{M}\mathrm{fld}}}
\newcommand{\Set}{{\mathcal{S}\mathrm{et}}}
\newcommand{\Cat}{{\mathcal{C}\mathrm{at}}}
\newcommand{\sSet}{{\mathrm{s}\mathcal{S}\mathrm{et}}}
\newcommand{\sDelta}{{\boldsymbol \Delta}}

\DeclareMathOperator{\Ev}{Ev}
\DeclareMathOperator{\Cv}{Cv}

\DeclareMathOperator{\Hom}{Hom}

\DeclareMathOperator{\Iso}{Iso}
\DeclareMathOperator{\Graph}{graph}
\DeclareMathOperator{\Id}{Id}
\DeclareMathOperator{\id}{id}
\DeclareMathOperator{\ass}{ass_G}
\DeclareMathOperator{\Ass}{ass}
\DeclareMathOperator{\Assc}{ass_{\mathrm{2-cat}}}
\DeclareMathOperator{\LUni}{lun}
\DeclareMathOperator{\RUni}{run}
\DeclareMathOperator{\luni}{lun_G}
\DeclareMathOperator{\runi}{run_G}

\DeclareMathOperator{\linv}{linv_G}
\DeclareMathOperator{\rinv}{rinv_G}

\DeclareMathOperator{\pr}{pr}

\DeclareMathOperator{\Kan}{Kan}
\DeclareMathOperator{\Nerve}{\mathcal{N}}

\DeclareMathOperator{\Comma}{\downarrow}
\DeclareMathOperator{\Empty}{{\_\!\_\,}}

\newcommand{\Lcat}[1]{\underrightarrow{M}}
\newcommand{\Lgpd}[1]{\underleftrightarrow{M}}

\newcommand{\Bibu}[1]{{\hat{#1}}} 
\newcommand{\wideBibu}[1]{{\widehat{#1}}} 

\renewcommand{\phi}{\varphi}

\newenvironment{xyc}{\gathered\xy}{\endxy\endgathered}

\theoremstyle{plain}
\newtheorem{Theorem}{Theorem}[section]
\newtheorem{Proposition}[Theorem]{Proposition}
\newtheorem{Corollary}[Theorem]{Corollary}
\newtheorem{Lemma}[Theorem]{Lemma}

\theoremstyle{definition}
\newtheorem{Definition}[Theorem]{Definition}

\theoremstyle{remark}

\newtheorem{Example}[Theorem]{Example}

\def\Copr(#1,#2){
  \save
  (0,0)+(#1,#2)*{}="Node"; 
  (-1,1)+(#1,#2)*{}="In1";
  (1,1)+(#1,#2)*{}="In2";
  (0,-1)+(#1,#2)*{}="Out";
  (0,1)+(#1,#2)*\cir(1,0){d^u};
  "Node"; "Out" ** \dir{-};
  \restore
}

\def\Prod(#1,#2){
  \save
  (0,0)+(#1,#2)*{}="Node"; 
  (-1,-1)+(#1,#2)*{}="Out1";
  (1,-1)+(#1,#2)*{}="Out2";
  (0,1)+(#1,#2)*{}="In";
  "In"; "Node" ** \dir{-};
  (0,-1)+(#1,#2)*\cir(1,0){u^d};
  \restore
}

\def\Iden(#1,#2){
  \save
  (0,0)+(#1,#2)*{}="Node"; 
  (0,1)+(#1,#2)*{}="In";
  "In"; "Node" ** \dir{-};
  \restore
}

\def\Term(#1,#2){
  \save
  {\ar@{{*}} (#1,#2)*{}; (#1,#2)*{}};
  \restore
}

\def\Coun(#1,#2){
  \save
  (0,0)+(#1,#2)*{}="Node"; 
  (0,1)+(#1,#2)*{}="In";
  {\ar@{{*}} "Node"; "In" };
  \restore
}

\def\Unit(#1,#2){
  \save
  (0,0)+(#1,#2)*{}="Node"; 
  (0,1)+(#1,#2)*{}="In";
  {\ar@{{*}} "In"; "Node"};
  \restore
}

\def\Anti(#1,#2){
  \save
  (0,0)+(#1,#2)*{}="Node"; 
  (0,0.5)+(#1,#2)*\cir<2pt>{}="Mid";
  (0,1)+(#1,#2)*{}="In";
  "In"; "Mid" ** \dir{-};
  "Mid"; "Node" ** \dir{-};
  \restore
}

\def\Coev(#1,#2){
  \save
  (0,1)+(#1,#2)*\cir(1,0){d^u};
  \restore
}

\def\Eval(#1,#2){
  \save
  (0,-1)+(#1,#2)*\cir(1,0){u^d};
  \restore
}

\def\CoprBig(#1,#2)#3{
}

\def\ProdBig(#1,#2)#3{
}

\def\IdenBig(#1,#2)#3{
}

\def\Flip(#1,#2){
  \save
  (1,-1)+(#1,#2)*{}; (-1,1)+(#1,#2)*{} 
  ** \crv{ (1,0)+(#1,#2) & (-1,0)+(#1,#2) };
  (-1,-1)+(#1,#2)*{}; (1,1)+(#1,#2)*{}  
  ** \crv{ (-1,0)+(#1,#2) & (1,0)+(#1,#2) };
  \restore
}

\def\Flip(#1,#2){
  \save
  @(
  @+ {(1,0)+(#1,#2)} @+ {(-1,0)+(#1,#2)};
  (1,-1)+(#1,#2)*{}; (-1,1)+(#1,#2)*{} ** \crvs{};
  (1,1)+(#1,#2)*{}; (-1,-1)+(#1,#2)*{} ** \crvs{};
  @i @)
  \restore
}

\def\Morp(#1,#2)#3{
  \save
  (0,2)+(#1,#2)*{}="Out";
  (0,1)+(#1,#2)*+<2pt>[F-]{\text{\small #3}}="Center";
  (0,0)+(#1,#2)*{}="In";
  "In"; "Center" ** \dir{-};
  "Out"; "Center" ** \dir{-};
  \restore
}

\def\MorpB(#1,#2)#3{
}

\def\MorpM(#1,#2){
  \save
  (-1.2,2.2)+(#1,#2)*{}="Out1";
  (1.2,2.2)+(#1,#2)*{}="Out2";
  (0,1)+(#1,#2)*+<3pt>[F-]{\text{\small M}}="Center";
  (0,-0.2)+(#1,#2)*{}="In";
  "In"; "Center" ** \dir{-};
  "Out1"; "Center" ** \dir{-};
  "Out2"; "Center" ** \dir{-};
  \restore
 }

\def\MorpN(#1,#2){
  \save
  (-0.8,-0.2)+(#1,#2)*{}="Out1";
  (0.8,-0.2)+(#1,#2)*{}="Out2";
  (0,1)+(#1,#2)*+<3pt>[F-]{\text{\small N}}="Center";
  (0,2.2)+(#1,#2)*{}="In";
  "In"; "Center" ** \dir{-};
  "Out1"; "Center" ** \dir{-};
  "Out2"; "Center" ** \dir{-};
  \restore
 }

\begin{document}
\title{Stacky Lie groups}

\author{Christian Blohmann}

\address{Department of Mathematics, University of California, Berkeley, CA 94720, USA.\vspace*{-2ex}}
\address{Fakult\"at f\"ur Mathematik, Universit\"at Regensburg, 93040 Regensburg, Germany.}

\email{blohmann@math.berkeley.edu}

\thanks{This work was supported by a Marie Curie Fellowship of the European Union}





\begin{abstract}
Presentations of smooth symmetry groups of differentiable stacks are studied within the framework of the weak 2-category of Lie groupoids, smooth principal bibundles, and smooth biequivariant maps. It is shown that principality of bibundles is a categorical property which is sufficient and necessary for the existence of products. Stacky Lie groups are defined as group objects in this weak 2-category. Introducing a graphic notation, it is shown that for every stacky Lie monoid there is a natural morphism, called the preinverse, which is a Morita equivalence if and only if the monoid is a stacky Lie group. As example we describe explicitly the stacky Lie group structure of the irrational Kronecker foliation of the torus.
\end{abstract}

\maketitle

\begin{center}
  Dedicated to the memory of my friend Des Sheiham
\end{center}

\section{Introduction}

While stacks are native to algebraic geometry \cite{DeligneMumford:irreducibility} \cite{Giraud:Cohomologie}, there has recently been increasing interest in stacks in differential geometry. In analogy to the description of moduli problems by algebraic stacks, differentiable stacks can be viewed as a smooth structures describing generalized quotients of smooth manifolds. In this way, smooth stacks are seen to arise in genuinely differential geometric problems. For example, proper \'etale differentiable stacks, the differential geometric analog of Deligne-Mumford stacks, are a conceptually powerful way to describe orbifolds \cite{MoerdijkPronk:Orbifolds} \cite{Moerdijk:Orbifolds}. Another example comes from the observation that while not every Lie algebroid can be integrated by a Lie groupoid \cite{CrainicFernandes:Integrability}, it can always be integrated by a smooth stack of Lie groupoids \cite{ZhuTsegn:Integrating} which shows that smooth stacks are the appropriate setting for Lie's third theorem for algebroids. A third class of examples is provided by the appearance of bundle gerbes with (generalized) connections in string theory \cite{Freed:Higher}, gerbes over orbifolds \cite{LupercioUribe:Gerbes}, and the application of gerbes to twisted $K$-theory \cite{TuXuLaurent-Gengoux:Twisted}, which suggested the definition of gerbes over general differentiable stacks \cite{BehrendXu:Differentiable} \cite{BehrendXu:S1-bundles} \cite{Heinloth:Notes} \cite{Metzler:Topological}.

A differentiable stack is a stack over the category of manifolds together with an atlas, which is a surjective representable submersion from a manifold to the stack. From an atlas we obtain a presentation of the stack by a Lie groupoid. Conversely, given a Lie groupoid we obtain the stack it presents as the category of torsors. It is common folklore that isomorphic differentiable stacks correspond to Morita equivalent Lie groupoids \cite{BehrendXu:Differentiable} \cite{Zhu:n-groupoids}. The rigorous statement is that there is a weak 2-equivalence between the 2-category of differentiable stacks, fiber-preserving functors, and natural transformations on the one hand and the weak 2-category (i.e., bicategory) of Lie groupoids, right principal bibundles, and smooth biequivariant maps on the other hand. 

Many aspects of stacks given by fibered categories are inherently non-constructive. For example, in general there is no constructive choice of pull-backs, so we need to use the axiom of choice for proper classes for a section of the fiber functor to exist. Therefore, concrete calculations have to pass through a presentation. Moreover, in virtually all applications differentiable stacks arise through their presentations. This is one of the reasons why the focus in this paper is on the presentations by Lie groupoids.

The motivation of this work is to understand on a differential geometric level the symmetries of differentiable stacks. Here is a guiding example: Consider the action of the integers on the plane by rotation, $\bbZ \times \bbC \rightarrow \bbC$, $(k, z) \mapsto k \cdot z := e^{i \lambda k}z$ for $\lambda \in [0,2\pi]$. When $\lambda/2\pi$ is rational, the quotient $\bbC/\bbZ$ by this action looks like the surface of an infinite cone, which has an obvious $S^1$-symmetry. But if $\lambda/2\pi$ is irrational each orbit wraps around densely in a circle. The resulting quotient space is nowhere differentiable, the topology highly non-Hausdorff, and it does not seem to have any smooth symmetry. This seemingly pathological situation becomes innocuous if we instead describe the quotient by the differentiable stack naturally presented by the action groupoid. As illustration of the machinery set up in this paper we will show that this stack does have a smooth symmetry given by a stacky Lie group, which is presented by the reduced groupoid of the Kronecker foliation of the 2-torus. It has been shown elswhere that this stacky group structure gives rise to a hopfish structure on the noncommutative torus algebra \cite{BlohmannTangWeinstein:Hopfish}. Another motivation for this work was to understand stacky Lie groups as the geometric counterparts of hopfish algebras \cite{TangWeinsteinZhu:Hopfish}. 

The content of this paper is as follows: Section~\ref{sec:Stacks} presents review material in order to make the paper reasonably self-contained and fix the terminology. We give a brief overview of differentiable stacks and how they are presented by Lie groupoids. We review the concept of smooth right principal bibundles of Lie groupoids and explain how they become the 1-morphisms of a weak 2-category. We recall how to associate bibundles to homomorphisms of groupoids. Finally, we explain how to get back a differentiable stack from a Lie groupoid. 

In Section~\ref{sec:Bibundles}, we give three alternative characterizations of principality of smooth bibundles. First, we show that right principality is equivalent to the existence of a groupoid-valued bibundle pairing (Proposition~\ref{th:freetransitive} and Corollary~\ref{th:submersion}), which can be viewed as geometric analog of the algebra-valued inner product of a Hilbert module. Second, we show that the structure of a bibundle can be expressed in terms of a category, which we call the linking category (Proposition~\ref{th:linking}). We give the conditions for the bibundle to be principal in terms of the linking category (Proposition~\ref{th:linkingprincipal}). If the bibundle is biprincipal we even get a linking groupoid (Proposition~\ref{th:linkinggroupoid}), which is the geometric analog of the linking algebra of a Hilbert bimodule. Third, we show that right principality is a categorical property, being equivalent to the 2-commutativity of two diagrams in the weak 2-category of groupoids, bibundles, and biequivariant maps (Theorem~\ref{th:rightprincipal}). Using this third characterization we show that right principality is a sufficient (Proposition~\ref{th:CartesianCat}) and necessary (Proposition~\ref{th:CartesianCat2}) condition for finite weak categorical products to exist. This is essential, because we need weak categorical products in order to have the notion of weak group objects, which we call stacky Lie groups. We further observe that the groupoid inverse can be viewed as a weak rigid structure on the weak monoidal 2-category of groupoids and bibundles (Proposition~\ref{th:rigid}). We introduce a graphic notation which allows for intuitive and concise proofs. 

In Section~\ref{sec:Groups}, we observe that using the evaluation and coevaluation bibundles of the rigid structure we can construct for every stacky Lie monoid a smooth bibundle, called the stacky preinverse, and show that a stacky Lie monoid is a stacky Lie group if and only if this preinverse is a Morita equivalence of Lie groupoids (Theorem~\ref{th:HopfishInverse}). The proof is given in the graphic notation which makes it purely categorical and which justifies rigorously why this property is the geometric analog of one of the axioms of hopfish algebras. Finally, we describe explicitly the stacky group structure of the Kronecker foliation of torus.

\section{Differentiable stacks and Lie groupoids}
\label{sec:Stacks}

A differentiable stack is a stack over the category of manifolds which has an atlas. We will briefly explain these notions. For a rigorous and detailed treatment see \cite{BehrendXu:Differentiable}.

\subsection{From differentiable stacks to Lie groupoids}

\begin{Definition}
Let $\catX$, $\catB$ be categories and $F: \catX \rightarrow \catB$ a functor. We say that a morphism $X \rightarrow Y$ in $\catX$ lies over a morphism $U \rightarrow V$ in $\catB$ if $F(X \rightarrow Y) = U \rightarrow V$. $F$ is called a \textbf{fibration} if the following holds:
\begin{itemize}

 \item[(F1)] For every morphism $U \rightarrow V$ in $\catB$ and every object $Y \in \catX$ over $V$ there is a morphism $X \rightarrow Y$ over $U \rightarrow V$. 

 \item[(F2)] Every morphisms $X \rightarrow Z$ over $U \rightarrow V \rightarrow W$ factorizes uniquely through every morphism $Y \rightarrow Z$ over $V \rightarrow W$.

\end{itemize}
\end{Definition}
We can express (F1) and (F2) by the following diagrams:
\begin{equation*}
\xymatrix@+1.5ex{\exists X \ar@{.>}[r]^{\exists} \ar@{|.>}[d] & Y \ar@{|->}[d] 
\\ U \ar[r] & V
}
\qquad\qquad\qquad
\xymatrix@-3ex{ 
X \ar@{.>}[rr]^{\exists !} \ar[dr] \ar@{|->}[dd] & & 
Y \ar[dl] \ar@{|->}[dd] \\ & Z  \ar@{|->}[dd] \\
U \ar'[r][rr] \ar[dr] & & V \ar[dl] \\ & W 
}
\end{equation*}
We also call $\catX$ a category fibered over $\catB$. The set $F_U := F^{-1}(1_U)$ of arrows over $1_U$, called the \textbf{fiber} over $U$, is a subcategory of $\catX$ which is a groupoid by (F2). For this reason $\calX$ is also called a category fibered in groupoids. 

Let now $\catB$ be equipped with a Grothendieck topology. Here we will only consider the category of manifolds as base category $\catB = \Mfld$ with the Grothendieck topology of coverings by smooth open embeddings. Given a covering family $U_i \rightarrow U$, we denote $U_{ij} := U_i \times_U U_j = U_i \cap U_j$ and $U_{ijk} := U_i \times_U U_j \times_U U_k$. 


Let $\alpha_i: X_i \rightarrow X$ and $\beta_i: Y_i \rightarrow Y$ be families of morphisms over a covering family $U_i \rightarrow U$. Given a morphism $\phi:X \rightarrow Y$ in $F_U$, there is by (F2) a unique family of morphisms $\phi|_{U_i}: X_i \rightarrow Y_i$ such that $\beta_i \circ \phi|_{U_i} = \phi \circ \alpha_i$, called the \textbf{pull-backs} of $\phi$ to $U_i$. Using families of morphisms $\alpha_{ij}: X_{ij} \rightarrow X_i$ and $\beta_{ij}: Y_{ij} \rightarrow Y_i$ over the induced covering $U_{ij} \rightarrow U_i$ we can iterate this construction to give pull-backs of pull-backs $\phi|_{U_{ij}} = (\phi|_{U_i})|_{U_{ij}}$ etc.

\begin{Definition}
Let $F: \catX \rightarrow \catB$ a fibration over a category $\catB$ with a Grothendieck topology. $\catX$ is called a \textbf{stack} if it satisfies the following gluing axioms:
\begin{itemize}
 \item[(S1)] Let $\phi,\psi: X \rightarrow Y$ be isomorphisms in $F_U$. If $\phi|_{U_i} = \psi|_{U_i}$ for all $i$, then $\phi = \psi$.

 \item[(S2)] Let $X$, $Y$ be objects over $U$ and $\phi_i : X_i \rightarrow Y_i$ be a family of isomorphisms. If $\phi_{i}|_{U_{ij}} = \phi_{j}|_{U_{ij}}$ for all $i$, $j$, then there is an isomorphism $\phi: X \rightarrow Y$ such that $\phi_i = \phi|_{U_i}$.

\item[(S3)] Let $\phi_{ij}: X_{ij} \rightarrow X_{ji}$ be a family of morphisms. If $\phi_{ij}|_{U_{ijk}} \circ \phi_{jk}|_{U_{ijk}} = \phi_{ik}|_{U_{ijk}}$ for all $i$, $j$, $k$, then there is an object $X$ over $U$ and morphisms $\phi_i : X_i \rightarrow X$ such that $\phi_i|_{U_{ij}} = \phi_{ij}$.
\end{itemize}
\end{Definition}

Axioms (S1) and (S2) are often expressed by saying that ``isomorphism is a sheaf'' and (S3) by saying that ``all descent data are effective''.

Every manifold $X$ gives rise to a stack $\Stack{X}$ over the category of manifolds given by the category of smooth maps to $X$. A morphism from a smooth map $U \rightarrow X$ to another smooth map $V \rightarrow X$ is a commutative triangle of smooth maps
\begin{equation*}
\xymatrix@-1ex@C-2ex{
  U \ar[rr]\ar[dr] && V \ar[dl] \\
 & X } 
\end{equation*}
The fiber functor maps such a triangle to the smooth map $U \rightarrow V$. A morphism of stacks $F: \catX \rightarrow \Mfld$ and $G: \catY \rightarrow \Mfld$ is a functor $f: \catX \rightarrow \catY$ which preserves fibers, $G \circ f = F$. A morphism of stacks $\Stack{X} \rightarrow \Stack{Y}$ is the same as a smooth map from $X$ to $Y$. Therefore, the category of manifolds is a full subcategory of the category of stacks over manifolds, so one can simply speak of $\Stack{X}$ as a manifold, of a morphism $\Stack{X} \rightarrow \Stack{Y}$  being a submersion, etc.

Natural transformations of fiber preserving functors of stacks can be viewed as 2-mor\-phisms. This gives rise to the strict 2-category of stacks, fiber preserving functors, and natural transformations. Stacks are called \textbf{isomorphic} if they are weakly isomorphic in this 2-category, in other words, if there is a fiber preserving equivalence of categories. A stack over manifolds is called \textbf{representable} if it is isomorphic to the stack of a manifold. 

Given two morphisms of stacks $l: \catX \rightarrow \catZ$ and $r:\catY \rightarrow \catZ$ we have the \textbf{2-fiber product} $\catX \times_{\catZ} \catY$ defined as terminal object making the diagram
\begin{equation*}
\xymatrix{
 \catX \times_{\catZ} \catY \ar[r] \ar[d] \drtwocell<\omit> & \catY \ar[d]^{r}\\
 \catX \ar[r]^{l} & \catZ}
\end{equation*}
\textbf{2-commutative}, which means that there is a 2-isomorphism from $\catX \times_{\catZ} \catY \rightarrow \catY \rightarrow \catZ$ to $\catX \times_{\catZ} \catY \rightarrow \catX \rightarrow \catZ$ indicated in the diagram by the diagonal double-arrow. Explicitly, an object of $\catX \times_{\catZ} \catY$ is a triple $(X,\alpha,Y)$ consisting of objects $X \in \catX$ and $Y \in \catY$ with $l(X) = r(Y)$, and an isomorphism $l(X) \stackrel{\alpha}{\leftarrow} r(Y)$ in $\catZ$. 

A morphism of stacks $\catX \rightarrow \catY$ is called \textbf{representable} if for every manifold $T$ and morphism of stacks $\Stack{T} \rightarrow \catY$ the 2-fiber product $\catX \times_{\catY} \Stack{T}$ is representable. A representable morphism is called a \textbf{submersion} if the induced morphism $\catX \times_{\catY} \Stack{T} \rightarrow \Stack{T}$, which is a smooth map of manifolds, is a submersion. A morphism of stacks $\catX \rightarrow \catY$ is an \textbf{epimorphism} if for every $\Stack{T} \rightarrow \catY$, where $T$ is a manifold, there is a surjective submersion $S \rightarrow T$ and a morphism of stacks $\Stack{S} \rightarrow \catX$ such that the diagram
\begin{equation*}
\xymatrix{
 \Stack{S} \ar[r] \ar[d] \drtwocell<\omit> & \Stack{T} \ar[d]\\
 \catX \ar[r] & \catY}
\end{equation*}
is 2-commutative. Putting all these properties together, a morphism $\catX \rightarrow \catY$ is called a \textbf{surjective representable submersion} if it is an epimorphism which is a representable submersion.

\begin{Definition}
  A stack $\catX \rightarrow \Mfld$ is called \textbf{differentiable} if there is a manifold $S$ and a surjective representable submersion $\Stack{S} \rightarrow \catX$, called an \textbf{atlas}.
\end{Definition}

Given an atlas $s:\Stack{S} \rightarrow \catX$, consider the 2-fiber product
\begin{equation*}
\xymatrix{
 \Stack{S} \times_{\catX} \Stack{S} \ar[r] \ar[d] \drtwocell<\omit> & \Stack{S} \ar[d]\\
 \Stack{S} \ar[r] & \catX} 
\end{equation*}
Explicitly, $\Stack{S} \times_{\catX} \Stack{S}$ is given by pairs of maps $\mu:U \rightarrow S$ and $\nu:U \rightarrow S$ together with an isomorphism $s(\mu) \stackrel{\phi}{\leftarrow} s(\nu)$ in $\catX$. Composition of isomorphisms equips the pull-back stack with a natural groupoid structure. It follows immediately from the defining properties of an atlas that the pull-back stack is representable and that the projections $\Stack{G_1} \cong \Stack{S} \times_{\catX} \Stack{S} \rightrightarrows \Stack{S} =: \Stack{G_0}$ are surjective submersions. The conclusion is that $G_1 \rightrightarrows G_0$ is a Lie groupoid, called a \textbf{presentation} of the differentiable stack.
\begin{Example}
  A representable stack $\Stack{X}$ is its own atlas. The Lie groupoid presentation is given by $X_1 = X_0 = X$.
\end{Example}

\begin{Example}
  Let $G$ be a Lie group. $BG$, the category of right principal $G$-bundles and smooth $G$-equivariant bundle maps, is a stack over manifolds. An atlas is given by a point, $s: \Stack{*} \rightarrow BG$, $s(U \rightarrow {*}) := U \times G$. The presentation is given by the Lie group $G$ viewed as Lie groupoid with one object.
\end{Example}

\subsection{Lie groupoids and bibundles}

A groupoid $G = (G_1 \rightrightarrows G_0)$ is a category in which each arrow in $G_1$ has an inverse. Because there seems to be no agreement on whether the source or the target should be written on the right or the left we will call them the left and right \textbf{moment maps} from $G_1$ to $G_0$ and denote them by $l_G$ and $r_G$, respectively. If no confusion can arise we will drop the subscript.  When defining a groupoid internal to the category of differentiable manifolds, we have to require one of the moment maps to be a submersion so that the pullback $G_1 \times_{G_0}^{l,r} G_1$, on which the groupoid multiplication is defined, is a differentiable manifold. Applying the groupoid inverse shows that the other moment map is a submersion, too. A \textbf{Lie groupoid} is a groupoid in the category of differentiable manifolds with submersive moment maps \cite{Pradines:Theorie}. 

\begin{Example}
In our guiding example we considered the action of $\bbZ$ on $\bbC$ given by $k \cdot z = e^{i\lambda k} z$. The Lie groupoid $\gX = (\gX_1 \rightrightarrows \gX_0)$ presenting the differentiable stack $\bbC /\!/ \bbZ$ of the quotient by the group action is given by the manifolds $\gX_1 = \bbZ \times \bbC$ and $\gX_0 = \bbC$, with moment maps $l_{\gX}(k,z) = e^{ik\lambda} z$, $r_{\gX}(k,z) = z$, and groupoid multiplication $(k_1,z_1)(k_2,z_k) = (k_1 + k_2, z_2)$ whenever $z_1 = e^{ik_2\lambda} z_2$. The identity bisection is given by $1_z = (0,z)$.
\end{Example}


\begin{Definition}
Let $G$ be a Lie groupoid and $M$ a manifold with a smooth map $l_M: M \rightarrow G_0$, called the \textbf{moment map}. A differentiable map 
\begin{equation*}
  \rho: G_1 \times_{G_0}^{r_G,l_M} M \longrightarrow M
  \,,\qquad
  (g,m) \longmapsto g \cdot m \,,
\end{equation*}
for which $g \cdot (g' \cdot m) = gg' \cdot m$, whenever defined, is called a \textbf{left action} of $G$ on $M$.
\end{Definition}

We call a bundle with a left action of $G$ a left $G$-bundle. There is an obvious definition of a right action.

\begin{Definition}
  Let $G$ and $H$ be Lie groupoids. A manifold $M$ with a left $G$-action and a right $H$-action which commute, i.e., $(g \cdot m) \cdot h = g \cdot (m \cdot h)$, whenever defined, is called a smooth $G$-$H$ \textbf{bibundle}.
\end{Definition}

Most of the material about bibundles of this section can be found in \cite{HilsumSkandalis:Morphismes} \cite{Landsman:Quantized} \cite{Landsman:Lie_groupoids} \cite{Mrcun:Stability} \cite{Mrcun:Functoriality}. We can view a $G$-$H$ bibundle as a presentation of a relation between the differentiable stacks presented by $G$ and $H$ \cite{BlohmannWeinstein:Group-like}. 


\begin{Definition}
  Let $G$, $H$ be Lie groupoids and $M$,$N$ be $G$-$H$ smooth bibundles. A \textbf{morphism of smooth bibundles} is a smooth map $\phi: M \rightarrow N$ which is biequivariant, i.e., $l_N(\phi(m)) = l_M(m)$, $r_N(\phi(m)) = r_M(m)$ for all $m \in M$ and $\phi(g \cdot m) = g \cdot \phi(m)$, $\phi(m \cdot h) = \phi(m) \cdot h$ for all $m \in M$, $g \in G$, $h \in H$ for which the actions are defined.
\end{Definition}

Let $G$, $H$, $K$ be Lie groupoids, $M$ a $G$-$H$ bibundle, and $N$ a $H$-$K$ bibundle. Viewing a bibundle as relation of stacks suggests defining the composition of bibundles as
\begin{equation}
\label{eq:Bibucomp}
  M \circ N := (M \times_{H_0}^{r_M,l_N} N)/H \,,
\end{equation}
where the quotient is with respect to the diagonal action $(m,n) \cdot h = (m \cdot h, h^{-1} \cdot n)$ for $(m,n) \in M \times_{H_0} N$, $h \in H$ whenever defined. By construction, this composition of two smooth bibundles is a topological $G$-$K$ bibundle, which can fail to be smooth for two reasons: First, the pull-back $M \times_{H_0}^{r_M,l_N} N$ does not have a natural manifold structure unless one of the moment maps is a submersion. Second, the quotient by the groupoid action can fail to be smooth. For the composition of differentiable bibundles to be differentiable we have to require further properties.

\begin{Definition}
\label{def:principal}
  Let $G$, $H$ Lie groupoid. A smooth $G$-$H$ bibundle $M$ is called \textbf{right principal} if the following three conditions hold: 
\begin{itemize}
 \item[(P1)] The left moment map of $M$ is a surjective submersion.
 \item[(P2)] The right $H$-action is free.
 \item[(P3)] The right $H$-action is transitive on the $l_M$-fibers.
\end{itemize}
\end{Definition}

Note that (P1) is a property of the left action and (P2) is a property of the right action. Property (P3) involves the left fibration and the right action. The analogous definition of a left principal bibundle is obvious. We will give two useful alternative characterizations of right principal bibundles in the next section.

For topological bibundles, property (P1) becomes the requirement that the left moment map is a continuous surjection which has a local section at every point of the base. For set theoretic groupoids the condition (P1) becomes
\begin{itemize}
 \item[(P1)] The left moment map of $M$ is surjective.
\end{itemize}
The wording of conditions (P2) and (P3) remains unchanged.

When composing two smooth principal bibundles $M$ and $N$, property $(P1)$ of the right bibundle $N$ ensures that the pull-back $M \times_{H_0} N$ is naturally a manifold. Properties (P2) and (P3) of the left bibundle $M$ imply that the diagonal action of $H$ is free and proper (see Corollary~\ref{th:proper}) and, hence, that the quotient is smooth.

\begin{Proposition}
 Let $G$, $H$, $K$ be Lie groupoids, $M$ a smooth $G$-$H$ bibundle, and $N$ a smooth $H$-$K$ bibundle. If $M$ and $N$ are right principal then the composition $M \circ N$ is smooth right principal.
\end{Proposition}

Given two homomorphisms of bibundles $\phi : M \rightarrow M'$ and $\psi: N \rightarrow N'$ we have the map of the cartesian products $\phi \times \psi: M \times N \rightarrow M' \times N'$. Because of the $H$-equivariance of $\phi$ and $\psi$ this map restricts to a map of the pull-backs over $H_0$ which then descends to a map of the quotients $\phi \circ \psi : M \circ N \rightarrow N' \circ N'$. This composition of homomorphisms of bibundles is a functor to the category of bibundles and biequivariant maps, i.e., $(\phi_1 \circ \psi_1)(\phi_2 \circ \psi_2) = \phi_1 \phi_2 \circ \psi_1 \psi_2$ whenever defined. (Because we use the symbol $\circ$ for the functor of bibundle composition, we omit it for the concatenation of bibundle morphisms.)

The composition of bibundles is associative only up to natural isomorphisms of bibundles. Let $M$ be a $G$-$H$ bibundle, $L$ an $H$-$H'$ bibundle, and $N$ an $H'$-$K$ bibundle. We have a natural diffeomorphism
\begin{equation*}
  (M \times_{H_0} N) \times_{H_0'} L 
  \stackrel{\cong}{\longrightarrow}
  M \times_{H_0} (N \times_{H_0'} L) \,,
\end{equation*}
by identifying both manifolds as submanifold $X \subset M \times N \times L$. The actions of $H$ and $H'$ commute, so we get natural isomorphisms
\begin{equation*}
  (X/H)/H' 
  \stackrel{\cong}{\longrightarrow}
  X/(H \times H') 
  \stackrel{\cong}{\longrightarrow}
  (X/H')/H \,.
\end{equation*}
If we combine these natural isomorphisms we get a natural isomorphism of bibundles
\begin{equation*}
  \Ass_{M,N,L}: 
  (M \circ N) \circ L 
  \stackrel{\cong}{\longrightarrow}
  M \circ (N \circ L) \,,
\end{equation*}
called the associator. It is elementary to check that the associator satisfies the pentagon relation:
\begin{equation}
\label{eq:pentagondiag}
\begin{gathered}
\xymatrix@C-8ex{
 & ((A \circ B) \circ C) \circ D
  \ar[dl]_{\Ass_{A,B,C} \times \id_D ~~}
  \ar[dr]^{~~\Ass_{A\circ B,C,D}}
\\
 (A \circ (B \circ C)) \circ D 
  \ar[d]_{\Ass_{A,B\circ C, D}}
& & (A \circ B) \circ (C \circ D)
  \ar[d]^{\Ass_{A,B,C\circ D}}
\\
 A \circ ((B \circ C) \circ D)
  \ar[rr]^{\id_A \times \Ass_{B,C,D}}
& & A \circ (B \circ (C \circ D))
}
\end{gathered}
\end{equation}
The Lie groupoid $G$ with left and right groupoid multiplication is itself a $G$-$G$ bibundle which we denote by $\Id_G$. We can compose $\Id_G$ from the right with every $G$-$H$ bibundle $M$. It is straight-forward to check that the smooth maps
\begin{xalignat*}{2}
  \LUni_M:
  \Id_G \circ M &\longrightarrow M\,,&
  M &\longrightarrow \Id_G \circ M \\
  [g,m] &\longmapsto g \cdot m  \,,&
  m &\longmapsto [l_M(m), m] \,,
\end{xalignat*}
are mutually inverse and biequivariant. The isomorphism of bibundles $\LUni_M$ is called the left unit constraint. There is an analogous right unit constraint ${}_M\!\RUni: M \circ \Id_H \rightarrow M$. Both unit constraints are natural isomorphisms of bibundles. Again, it is elementary to check that the unit constraints satisfy the coherence relation:
\begin{equation*}
 \xymatrix@C-4ex{
(M \circ \Id_H) \circ N \ar[rr]^{\Ass_{MHN}} \ar[dr]_{{}_M\!\LUni}
 & & M \circ (\Id_H \circ N) \ar[dl]^{\RUni_N} \\
& M \circ N
}\quad.
\end{equation*} 
The conclusion is:

\begin{Proposition}
  Lie groupoids, smooth right principal bibundles, and smooth biequivariant maps of bibundles form a weak 2-category denoted by $\LieGrpdPrBibu$.
\end{Proposition}

The abbreviation $\LieGrpdPrBibu$ stands for ``stacky manifolds'', which is an intuitive way to think of this category. Throughout this article ``weak 2-category'' is a synonym for ``bicategory''. Occasionally, we denote the collection of all 1-morphisms from $G$ to $H$ by $\Hom_1(G,H)$ and the set of all 2-morphisms from $M$ to $N$ by $\Hom_2(M,N)$. A $G$-$H$ bibundle $M$ is called a \textbf{weak 1-isomorphism} if there is a $H$-$G$ bibundle $N$ such that $M \circ N \cong \Id_G$ and $N \circ M \cong \Id_H$ as bibundles. Two groupoids $G$ and $H$ are \textbf{weakly isomorphic} if there is a $G$-$H$ bibundle $M$ which is a weak isomorphism. This notion of weak isomorphism of groupoids is usually called \textbf{Morita equivalence}.

\subsection{Bundlization}


To every smooth homomorphism of Lie groupoids $\phi: G \rightarrow H$ we can associate a smooth $G$-$H$ bibundle $\Bibu{\phi}$ defined as the manifold 
\begin{equation*}
  \Bibu{\phi} := G_0 \times_{H_0}^{\phi_0,l_H} H_1
  = \{(x,h) \in G_0 \times H_1 \,|\, \phi(x) = l_H (h)\} \,,
\end{equation*}
with moment maps $l_{\Bibu{\phi}}(x,h) = x$, $r_{\Bibu{\phi}}(x,h) = r_H(h)$ and left and right groupoid actions
\begin{equation*}
  g\cdot (x,h) = (l_G(g), f(g)h) \,,\qquad (x,h) \cdot h' = (x,hh') \,,
\end{equation*}
defined whenever $r_G(g) = l_{\Bibu{\phi}}(x,h) = x$ and $l_H(h') = r_{\Bibu{\phi}}(x,h) = r_H(h)$. We will call $\Bibu{\phi}$ the \textbf{bundlization} of $\phi$.

The map $G_0 \rightarrow \Bibu{\phi}$, $x \mapsto (x, 1_{\phi_0(x)})$ is a smooth section of the left bundle $l_{\Bibu{\phi}}: \Bibu{\phi} \rightarrow G_0$. This shows that $l_{\Bibu{\phi}}$ is a surjective submersion. Moreover,  $h^{-1}h' \in H$ is the unique element which takes by the right $H$-action $(x,h) \in \Bibu{\phi}$ to $(x,h')$. The conclusion is:

\begin{Proposition}
The $G$-$H$ bibundle $\Bibu{\phi}$ is smooth right principal.
\end{Proposition}

It is straight-forward to check that bundlization $\phi \mapsto \Bibu{\phi}$ is a functor from the category of Lie groupoids and smooth homomorphisms to $\LieGrpdPrBibu$. This functor is not essentially full (see Proposition~\ref{th:section}).

Let $X$ be a smooth manifold. Then $X$ can be equipped with a trivial groupoid structure, $X_1 = X_0 := X$, $l_X = r_X := \id_X$. The bundlization of a smooth map $\phi: X \rightarrow Y$ is its graph $\Bibu{\phi} = \Graph(\phi) \subset X \times Y$ with the trivial groupoid action on both sides. The composition of graphs viewed as bibundles is the graph of the composition of the smooth maps. Conversely, for every smooth $X$-$Y$ bibundle $M$ right principality implies that the left moment map is a diffeomorphism. So we get a unique smooth map $\phi_M := l_M^{-1}r_M: X \rightarrow Y$ the bundlization of which is $M$. We conclude the following:

\begin{Proposition}
  The category of smooth manifolds is a full subcategory of $\LieGrpdPrBibu$.
\end{Proposition}

\subsection{From Lie groupoids to differentiable stacks}
\label{sec:GrpdsToStacks}

We now review briefly how to get back from a Lie groupoid to a differentiable stack \cite{BehrendXu:Differentiable}.

\begin{Definition}
  Let $G$ be a Lie groupoid. A right $G$-torsor is a 1-morphism in $\LieGrpdPrBibu$ from a manifold $X$ to $G$. Denote the collection of right $G$-torsors by
\begin{equation*}
  [G] := \{ S \in \Hom_1(X,G) \,|\, X \text{ is a manifold} \} \,.
\end{equation*}
\end{Definition}

In analogy to groups this is often called the classifying space of the Lie groupoid and denoted by $BG$. We will not use this terminology here to avoid confusion with delooping. First, we will indicate why $[G]$ is a differentiable stack.

$[G]$ is a category with smooth right $G$-equivariant bundle maps as morphisms. There is a functor $F$ which maps a smooth $G$-equivariant bundle map from $S \in \Hom_1(X,G)$ to $T \in \Hom_1(X,G)$ to the induced map $X \rightarrow Y$ of the base manifolds.

\begin{Proposition}
  Let $G$ be a Lie groupoid. Then $F: [G] \rightarrow \Mfld$ is a stack.
\end{Proposition}
\begin{proof}
Let $\phi: X \rightarrow Y$ be a morphism of smooth manifolds and let $S \in [G]$ be in the $Y$-fiber of $F$, i.e., $S$ is a right principal $Y$-$G$ bibundle. By bundlization of $\phi$ we get a right principal $X$-$Y$ bibundle $\Bibu{\phi}$ which we can compose with $S$ to obtain the bibundle $S|_X := \Bibu{\phi} \circ S = X \times_{Y}^{\phi,l_S} S$ which lies in the $X$-fiber of $F$. Moreover, the induced map $S|_X \rightarrow S$, $(x,s) \mapsto s$ is a $G$-equivariant bundle map which covers $\phi$. We conclude that axiom (F1) of a fibration holds.

Consider a $G$-equivariant bundle map $f_\phi: S|_X \rightarrow S$ covering the map $\phi: X \rightarrow Y$ and another such map $f_\psi: S|_Z \rightarrow S$ covering $\psi: Z \rightarrow Y$. If there is a map $\chi: X \rightarrow Z$ such that $\psi \circ \chi = \phi$, then the map $f_\chi: S|_X \rightarrow S|_Y$, $(x,s) \mapsto (\chi(x),s)$ is the unique $G$-equivariant bundle map covering $\chi$ such that $f_\psi \circ f_\chi = f_\phi$. We conclude that axiom (F2) of a fibration holds, as well.

We now turn to the gluing properties of $[G]$, to show that it is actually a stack. Let $\{U_i \hookrightarrow X \}$ be a covering by smooth open embeddings in the manifold $X$ and $U_{ij} := U_i \cap U_j$. Then $S|_{U_i}$ is the restriction of the bundle $l_S : S \rightarrow X$ to $U_i$. Given $S, T$ in the $X$-fiber of $F$, isomorphisms from $S$ to $T$ are equivariant diffeomorphisms. Denote by $\Iso(S,T)(U_i) := \Iso_2(S|_{U_i}, T|_{U_i})$ the isomorphisms between the restricted bundles. It is rather obvious that $U_i \mapsto \Iso(S, T)(U_i)$ is a sheaf on the Grothendieck topology of open smooth coverings. 

To show this explicitly, let $f,g: S \rightarrow T$ be equivariant diffeomorphisms of $S$, $T$ in the $X$-fiber of $[G]$. Denote by $f|_{U_i}$ the restriction of $f$ to $S|_{U_i}$. If all restrictions of $f$ and $g$ are equal, $f|_{U_i} = g|_{U_i}$ for all $i$, then $f = g$. This shows that axiom (S1) of a stack is satisfied. Moreover, if we have equivariant bundle isomorphisms $f_i : S|_{U_i} \rightarrow T|_{U_i}$ which are equal on the overlaps $f_i|_{U_{ij}} = f_j|_{U_{ij}}$ this defines a unique smooth isomorphism $f: S \rightarrow T$ which restricts to $f_i = f|_{U_i}$. This shows that axiom (S2) holds.

Next, we have to show that all descent data are effective. This is analogous to the construction of a principal fiber bundle from its locally trivializations. Given right principal $U_i$-$G$ bundles we first have to choose a refinement $\{U'_i \rightarrow X \}$ of the covering such that the restricted $U'_i$-$G$ bibundles $M_i$ have a section $\sigma_i$ of the left moment map $l_{M_i}: M_i \rightarrow U'_i$. This induces a smooth map $\phi_i: U'_i \rightarrow G_0$, $\phi_i(u) := r_{M_i}(\sigma(u))$ such that $M_i \cong U_i' \times_{G_0}^{\phi_i,l_G} G =: G|_{U_i'}$ as right principal $G$-bundles. Now we use the gluing properties of the restrictions $G|_{U_i'}$. Assume that we have equivariant diffeomorphisms $f_{ij} : G|_{U'_{ij}} \rightarrow G|_{U'_{ij}}$ whose restrictions to $U'_{ijk} := U'_i \cap U'_j \cap U'_k$ satisfy the group 2-cocycle condition $f_{jk} \circ f_{ij} = f_{ik}$. Then we have a pull-back bibundle $M$ and isomorphisms $f_i: M|_{U'_i} \rightarrow G|_{U'_i} \cong M_i$ such that $f_{ij} = f_j \circ f_i^{-1}$ on $U'_{ij}$. We conclude that axiom (S3) of a stack is satisfied.
\end{proof}

\begin{Proposition}
  Let $G$ be a Lie groupoid. Then $\Stack{G_0}$ is an atlas of $[G]$.
\end{Proposition}

\begin{proof}
Consider the morphism of fibered categories $\alpha: \Stack{G_0} \rightarrow [G]$ defined as follows: An object in $\Stack{G_0}$, given by a smooth function $\mu: U \rightarrow G_0$, is mapped to the bundlization $\alpha(\mu) := \Bibu{\mu} = U \times_{G_0}^{\mu,l_G} G_1$. A morphism in the category $\Stack{G_0}$ from $\mu: U \rightarrow G_0$ to $\nu: V \rightarrow G_0$ is given by a smooth map $\phi: U \rightarrow V$ such that $\mu \circ \phi = \nu$. This morphism is mapped to the $G$-equivariant map $\alpha(\phi) : \Bibu{\mu} \rightarrow \Bibu{\nu}$, $(u,g) \mapsto (\phi(u),g)$.

First, we show that $\alpha$ is an epimorphism. Let $T$ be a manifold and  $\beta: \Stack{T} \rightarrow [G]$ a functor of fibered categories. Let $T_i \rightarrow T$ be a covering family sufficiently fine such that all $\beta(T_i \rightarrow T)$ are trivial, i.e., there are maps $\sigma_i: T_i \rightarrow G_0$ such that $\beta(T_i \rightarrow T) \cong \Bibu{\sigma}_i$ as right $G$-torsors. Then $\Stack{T_i} \rightarrow \Stack{T} \rightarrow [G]$ maps a function $\mu: U \rightarrow T_i$ to $\wideBibu{\sigma_i \circ \mu}$. If we define $\tau_i: \Stack{T_i} \rightarrow \Stack{G_0}$ for all $i$ by $\tau_i(\mu) :=  \sigma_i \circ \mu$ we obtain 2-commutative squares
\begin{equation*}
\begin{gathered}
\xymatrix{
 \Stack{T_i} \ar[r] \ar[d]_{\tau_i} \drtwocell<\omit> & \Stack{T} \ar[d]^{\beta}\\
 \Stack{G_0} \ar[r]^{\alpha} & [G]}
\end{gathered}
\end{equation*}
Now, consider the 2-fiber product
\begin{equation*}
\begin{gathered}
\xymatrix{
 \Stack{G_0} \times_{[G]} \Stack{G_0} \ar[r] \ar[d] \drtwocell<\omit> 
 & \Stack{G_0} \ar[d]^{\alpha}\\
 \Stack{G_0} \ar[r]^{\alpha} & [G]}
\end{gathered}
\end{equation*}
It has been shown in \cite{BehrendXu:Differentiable} that if for an epimorphism $\alpha$ the 2-fiber product $\Stack{G_0} \times_{[G]}^{\alpha,\alpha} \Stack{G_0}$ is representable and if the projections on $\Stack{G_0}$ are submersions, then $\alpha$ is a representable submersion. We will show that the stack $\Stack{G_0} \times_{[G]} \Stack{G_0}$ is isomorphic to $\Stack{G_1}$.

The objects of $\Stack{G_0} \times_{[G]} \Stack{G_0}$ are given by a triple $(\mu,\phi,\nu)$ where $\mu,\nu: U \rightarrow G_0$ are smooth maps and $\Bibu{\mu} \stackrel{\phi}{\leftarrow} \Bibu{\nu}$ is an isomorphism of right $G$-torsors, given explicitly by $\phi(u,g) = (\phi_U(u,g),\phi_G(u,g))$. Right $G$-equivariance of $\phi$ reads
\begin{equation*}
\begin{split}
  \bigl( \phi_U(u,g),\phi_G(u,g) \bigr) 
  &= \phi(u,g) =  \phi(u,1_{\nu(u)}g) = \phi(u,1_{\nu(u)})g \\
  &= \bigl( \phi_U(u,1_{\nu(u)}), \phi_G(u,1_{\nu(u)})g \bigr) \,.
\end{split}
\end{equation*} 
We conclude that $\phi_U$ descends to a smooth map on $U$, also denoted by $\phi_U$. Because $\phi$ is a smooth isomorphism of $G$-torsors, $\phi_U$ is a diffeomorphism. For the $G$-component of $\phi$ we conclude that $r_G( \phi_G(u,1_{\nu(u)})) = \nu(u)$.

Given $(\mu,\phi,\nu) \in \Stack{G_0} \times_{[G]} \Stack{G_0}$, define the smooth map
\begin{equation*} 
  \eta_{(\mu,\phi,\nu)}: U \longrightarrow G \,,\qquad
  u \longmapsto \phi_G(u, 1_{\nu(u)}) \,.
\end{equation*}
Conversely, given a smooth map $\eta: U \rightarrow G$, define the triple of smooth maps
\begin{equation*}
\begin{aligned}
  U &\stackrel{\mu_\eta}{\longrightarrow} G_0 \\ 
  u &\longmapsto l_G(\eta(u)) 
\end{aligned}
\qquad
\begin{aligned}
  U &\stackrel{\nu_\eta}{\longrightarrow} G_0 \\
  u &\longmapsto r_G(\eta(u))
\end{aligned}
\qquad
\begin{aligned}
  \wideBibu{\nu_\eta} &\stackrel{\phi_\eta}{\longrightarrow} \wideBibu{\mu_\eta} \\
  (u,g) &\longmapsto (u, \eta(u)g) 
\end{aligned}
\end{equation*}
Because $\phi_\eta(u,gg') = (u,\eta(u)gg') = (u,\eta(u)g) \cdot g'$, the last map is a morphism of right $G$-torsors. The map $\wideBibu{\mu_\eta} \rightarrow \wideBibu{\nu_\eta}$, $(u,g) \mapsto (u,\eta(u)^{-1}g)$ is its inverse, so $\phi_\eta$ is an isomorphism.

So far we have constructed maps
\begin{equation}
\label{eq:StackIso}
\begin{aligned} 
  \Stack{G_0} \times_{[G]} \Stack{G_0} &\longrightarrow \Stack{G_1} \\\
  (\mu,\phi,\nu) &\longmapsto \eta_{(\mu,\phi,\nu)}
\end{aligned}
\qquad\qquad
\begin{aligned}
  \Stack{G_1} &\longrightarrow \Stack{G_0} \times_{[G]} \Stack{G_0} \\ 
  \eta &\longmapsto (\mu_\eta, \phi_\eta, \nu_\eta)
\end{aligned}\quad.
\end{equation}
It is straight-forward to check that these maps extend to morphisms of stacks. Composing the morphisms yields 
\begin{equation*}
\begin{aligned}
  \mu_{\eta_{(\mu,\phi,\nu)}}(u) 
  &= l_G\bigl(\eta_{(\mu,\phi,\nu)}(u) \bigr)
  = l_G\bigl( \phi_G(u, 1_{\nu(u)}) \bigr)
  = \mu(\phi_U(u)) \,,\\
  \nu_{\eta_{(\mu,\phi,\nu)}}(u) 
  &= r_G\bigl(\eta_{(\mu,\phi,\nu)}(u) \bigr)
  = r_G\bigl( \phi_G(u, 1_{\nu(u)}) \bigr)
  = \nu(u) \,,\\
  \phi_{\eta_{(\mu,\phi,\nu)}}(u,g) 
  &= \bigl(u,\eta_{(\mu,\phi,\nu)}(u) g\bigr)
  = \bigl(u, \phi_G(u, 1_{\nu(u)})g \bigr)
  = \bigl(u, \phi_G(u, g) \bigr) \,.
\end{aligned}
\end{equation*}
We conclude that $\phi_U$ is an isomorphism in $\Stack{G_0}$ from $\mu_{\eta_{(\mu,\phi,\nu)}}$ to $\mu$ which induces an isomorphism in $\Stack{G_0} \times_{[G]} \Stack{G_0}$ from $\phi_{\eta_{(\mu,\phi,\nu)}}$ to $\phi$. Hence, the morphisms of stacks of Eq.~\eqref{eq:StackIso} are an equivalence of fibered categories and, therefore, define an isomorphism of stacks $\Stack{G_0} \times_{[G]} \Stack{G_0} \cong \Stack{G_1}$.

At last, we observe that the induced projections $\Stack{G_1} \rightrightarrows \Stack{G_0}$ are given by the left and right moment maps of the Lie groupoid $G_1$, which are submersions. By Proposition~2.15 in \cite{BehrendXu:Differentiable} it follows that $\alpha$ is a representable submersion.
\end{proof}

We now briefly indicate how $G \mapsto [G]$ extends to a weak 2-functor from $\LieGrpdPrBibu$ to the 2-category of differentiable stacks, functors of fibered categories, and natural transformations.

Let $M$ be a smooth right principal $G$-$H$ bibundle. We have to show how $M$ induces a morphism of stacks, that is, a functor $[M]$ of fibered categories from $F_{[G]}: [G] \rightarrow \Mfld$ to $F_{[H]}: [H] \rightarrow \Mfld$. On an object $S \in [G]$ this functor is defined by $[M](S) := S \circ M$, which is compatible with the fibrations, $F_{[H]}([M](S)) = F_{[H]}(S \circ M) = F_{[G]}(S)$. A morphism in the category $[G]$ is a smooth equivariant map $f: S \rightarrow T$. On such a morphism the functor is defined by $[M](f) := (f \circ \id_M: S \circ M \rightarrow T \circ M)$. 

Let now $\phi: M \rightarrow N$ be a smooth equivariant map of $G$-$H$ bibundles, in other words a 2-morphism. This induces a natural transformation from the functor $[M]$ to $[N]$, by $[\phi]_S := \id_S \circ \phi: S \circ M \rightarrow S \circ N$, since for all smooth equivariant maps $f: S \rightarrow T$ the following diagram commutes:
\begin{equation}
\begin{gathered}
\xymatrix@+1ex{
S \circ M 
  \ar[r]^-{f \circ \id_M} 
  \ar[d]_{\id_S \circ \phi}
& T \circ M
  \ar[d]^{\id_T \circ \phi} 
\\ 
S \circ N 
\ar[r]^-{f \circ \id_N}
& T \circ N
}
\end{gathered}
\end{equation}
We have constructed a map of (weak) 2-categories
\begin{equation*}
\xymatrix@+1ex{G \rtwocell<5>^M_N{\phi} & H }
\quad \longmapsto \quad
\xymatrix@+1ex{[G] \rtwocell<5>^{[M]}_{[N]}{~[\phi]} & [H] }
\end{equation*}
Because $S \circ (M \circ N) \cong (S \circ M) \circ N$ and $f \circ \id_{(M \circ N)} = f \circ (\id_M \circ \id_N) = (f \circ \id_M) \circ \id_N$ the map $M \mapsto [M]$ is 1-functorial up to a natural transformation. And because $\id_S \circ (\phi\psi) = (\id_S \id_S) \circ (\phi\psi) = (\id_S \circ \phi)(\id_S \circ \psi)$ the map $\phi \mapsto [\phi]$ is 2-functorial. It is not difficult to show that this 2-functor is an equivalence:

\begin{Theorem}
  The weak 2-category $\LieGrpdPrBibu$ is 2-equivalent to the 2-category of differentiable stacks.
\end{Theorem}

\section{Categorical Properties of Lie groupoid bibundles}
\label{sec:Bibundles}

\subsection{The geometric analog of a Hilbert bimodule}

The following definition can be viewed as geometric analogue of the right algebra-valued inner product of a Hilbert bimodule. For brevity we call a bibundle \textbf{left submersive} if the left moment map is a submersion.

\begin{Definition}
 Let $G$, $H$ be Lie groupoids and $M$ a smooth, left submersive $G$-$H$ bibundle. An $H$-valued \textbf{bibundle pairing} is a smooth map
\begin{equation*}
  M \times_{G_0} M \longrightarrow H
  \,,\qquad
  (m,m') \longmapsto \langle m , m' \rangle \,,
\end{equation*}
which satisfies
\begin{itemize}
\item[(H1)] $\langle m , m' \cdot h \rangle = \langle m , m' \rangle h$
\item[(H2)] $\langle m , m' \rangle = \langle m' , m \rangle^{-1}$
\item[(H3)] $\langle m, m \rangle = 1_{r(m)}$
  \text{ and }
$\langle m , m' \rangle = \langle m , m'' \rangle$ implies $m' = m''$ 
\item[(H4)] $\langle g \cdot m , m' \rangle = \langle m , g^{-1} \cdot m' \rangle$
\end{itemize}
for all $m,m' \in M$, $g \in G$, $h \in H$ for which the expressions are defined.
\end{Definition}

Note, that we have to assume that $M$ is left submersive for $M \times_{G_0} M$ to be a manifold. For Hilbert bimodules we would have two $C^*$-algebras $A$ and $B$ instead of the groupoids $G$ and $H$, an $A$-$B$ bimodule $X$ instead of a bibundle $M$, and a $B$-valued inner product $\langle \xi, \xi' \rangle_B \in B$ for all $\xi,\xi' \in X$ instead of an $H$-valued pairing. Property (H1) is analogous to the Hilbert inner product being $B$-linear in the second argument, $\langle \xi, \xi' \cdot b \rangle_B = \langle \xi, \xi' \rangle_B b$. Property (H2) is analogous to the inner product being hermitian, $\langle \xi, \xi' \rangle_B = \langle \xi', \xi \rangle_B^*$. Property (H3) is analogous to the inner product being non-degenerate. And property (H4) is analogous to the action of $A$ being a $C^*$-action by adjointable operators, $\langle a \cdot \xi, \xi' \rangle_B = \langle \xi, a^* \cdot \xi' \rangle_B$.

\begin{Example}
  Viewing a Lie groupoid $G$ as $G$-$G$ bibundle, the pairing is given by $\langle g, g' \rangle = g^{-1}g'$.
\end{Example}

\begin{Example}
  For the bundlization $\Bibu{\phi}$ of a homomorphism $\phi:G \rightarrow H$ of Lie groupoids the pairing is given by $\langle (x,h), (x',h') \rangle = h^{-1}h'$.
\end{Example}

\begin{Proposition}
\label{th:freetransitive}
  Let $G$, $H$ be Lie groupoids and $M$ a smooth, left submersive $G$-$H$ bibundle. The following are equivalent:
\begin{itemize}
 \item[(i)] $M$ has an $H$-valued bibundle pairing.
 \item[(ii)] $M$ satisfies (P2) and (P3), i.e., the $H$-action is free and transitive on left fibers.
\end{itemize}
If these conditions hold, the bibundle pairing is unique.
\end{Proposition}

\begin{proof}
(i) $\Rightarrow$ (ii): Assume that for some $m \in M$, $h \in H$ we have $m \cdot h = m$. Then by (H1) we have $\langle m, m \rangle = \langle m, m \cdot h \rangle = \langle m, m \rangle h$ from which we conclude that $h = 1_{r(m)}$, so $H$ acts freely. Properties (H1), (H2), and (H3) imply that $\langle m' , m \cdot \langle m, m' \rangle \rangle = \langle m' , m \rangle \langle m, m' \rangle = 1_{r(m')} = \langle m', m' \rangle$. The second part of (H3) then implies that $m \cdot \langle m, m' \rangle = m'$, so $H$ acts transitively.

(ii) $\Rightarrow$ (i): We now assume that the action of $H$ is free and transitive on left fibers. For any two elements $m,m' \in M$ in the same left fiber $l(m) = l(m')$ we then have a unique element $\langle m , m' \rangle \in H$ such that $m \cdot \langle m , m' \rangle = m'$. Then $(m,m') \mapsto \langle m , m' \rangle$ defines a map from $M \times_{G_0} M \rightarrow H$, which satisfies properties (H1) through (H4). But we still have to show that the pairing is differentiable. For this we observe that the map 
\begin{equation*}
  \phi: M \times_{H_0} H \longrightarrow M \times_{G_0} M
  \,, \qquad
 (m,h) \longmapsto (m,m\cdot h)
\end{equation*}
is a local diffeomorphism because the $H$-action is free and transitive. The map
\begin{equation}
\label{eq:psiproper}
  \psi: M \times_{G_0} M \longrightarrow M \times_{H_0} H
  \,,\qquad
  (m,m') \mapsto (m, \langle m, m' \rangle )
\end{equation}
is the inverse of $\phi$. Because $\phi$ is a local diffeomorphism its inverse is differentiable. The bibundle pairing is $\psi$ followed by the projection to $H$, so it is also differentiable.

Let $\langle\!\langle~,~\rangle\!\rangle$ be another $H$-valued bibundle pairing on $M$. Then
\begin{equation*}
 \langle\!\langle m, m'\rangle\!\rangle
 = \langle\!\langle m, m \cdot \langle m, m' \rangle \rangle\!\rangle
 = \langle\!\langle m, m \rangle\!\rangle \langle m, m' \rangle 
 = \langle m, m' \rangle \,,
\end{equation*}
which shows that the bibundle pairing is unique.
\end{proof}

The existence of a bibundle pairing does not imply that the left moment map of the bibundle is surjective. For this we need an additional condition:

\begin{Corollary}
\label{th:submersion}
Let $G$, $H$ be Lie groupoids and $M$ a smooth, left submersive $G$-$H$ bibundle. Then $M$ is right principal if and only if it has an $H$-valued bibundle pairing and if the left moment map of $M$ has a smooth local section at every point $x \in G_0$.
\end{Corollary}
\begin{proof}
  A smooth map has local sections everywhere if and only if it is a surjective submersion.
\end{proof}

\begin{Corollary}
\label{th:proper}
  Let $M$ be a smooth, left submersive $G$-$H$ bibundle which satisfies (P2) and (P3). Then the right $H$-action is proper.
\end{Corollary}
\begin{proof}
 The map $\psi$ constructed in Eq.~\eqref{eq:psiproper} is a diffeomorphism and thus proper.
\end{proof}

This corollary ensures that the composition of smooth principal bibundles is smooth. Using the concept of bibundle pairings we can give a concise proof of the following fact:

\begin{Proposition}
\label{th:section}
 Let $G$, $H$ be Lie groupoids and $M$ a smooth right principal $G$-$H$ bibundle. $M$ is isomorphic to the bundlization of a smooth homomorphism of Lie groupoids if and only if there is a smooth section of the left moment map of $M$.
\end{Proposition}
\begin{proof}
 If $M = \Bibu{\phi} = G_0 \times^{\phi_0,l_H}_{H_0} H_1$ is the bundlization of a smooth homomorphism of Lie groupoids $\phi: G \rightarrow H$, then $\sigma: G_0 \rightarrow M$, $\sigma(x) = (x,1_{\phi_0(x)})$ is a smooth section. 

Conversely, let $\sigma$ be a smooth section of $l_M$. Consider the smooth map
\begin{equation*}
  \phi: G \longrightarrow H
  \,,\qquad
  g \longmapsto \bigl\langle \sigma(l(g)), g \cdot \sigma(r(g)) \bigr\rangle \,.
\end{equation*}
A direct calculation,
\begin{equation*}
\begin{split}
  \phi(g_1 g_2) 
  &= \bigl\langle \sigma(l(g_1 g_2)), g_1 g_2 \cdot \sigma(r(g_1 g_2)) \bigr\rangle \\
  &= \bigl\langle g_1^{-1} \cdot \sigma(l(g_1)), g_2 \cdot \sigma(r(g_2)) \bigr\rangle \\
  &= \bigl\langle g_1^{-1} \cdot \sigma(l(g_1)), l(g_2) \cdot 
     \langle \sigma(l(g_2)), g_2 \cdot \sigma(r(g_2)) \rangle \bigr\rangle \\
  &= \bigl\langle g_1^{-1} \cdot \sigma(l(g_1)), r(g_1) \bigr\rangle
     \bigl\langle \sigma(l(g_2)), g_2 \cdot \sigma(r(g_2)) \bigr\rangle\\
  &= \phi(g_1) \phi(g_2) \,,
\end{split}
\end{equation*}
shows that $\phi$ is a homomorphism of Lie groupoids. For $x \in G_0$ we have $\phi_0(x) = r_H(\sigma(x))$, so the bundlization is given by
\begin{equation*}
 \Bibu{\phi} = \{ (x,h) \in G_0 \times H_1 \,|\, r_H(\sigma(x)) = l_H(h) \} \,.
\end{equation*}
The map $\Bibu{\phi} \rightarrow M$, $(x,h) \mapsto \sigma(x) \cdot h$ is a smooth biequivariant map with inverse $M \rightarrow \Bibu{\phi}$, $m \mapsto (l_M(m), \langle \sigma(l_M(m)), m \rangle )$. 
\end{proof}

This distinguishes Lie groupoids or topological groupoids from set-the\-o\-re\-tic or categorical groupoids, where every surjective map has a section (assuming the axiom of choice). Observe that $l_M : M \rightarrow G_0$ is only a surjective submersion and in general not a fiber bundle. Therefore, it can fail to have a smooth section even if $G_0$ is contractible.

\subsection{The linking category and the linking groupoid}

In this section we consider only set-theoretic groupoids and bibundles without a differentiable structure. Intuitively, a $G$-$H$ groupoid bibundle $M$ is a collection of arrows from $G_0$ to $H_0$ such that the collection of all arrows in $G_1$, $M$, $H_1$ are a category. Let us make this statement more precise. 

Let $G$, $H$ be groupoids and $M$ a set with two maps $l_M: M \rightarrow G_0$ and $r_M: M \rightarrow H_0$. We can ask whether we can extend the groupoids $G$ and $H$ to a category with arrows and objects 
\begin{equation*}
  \Lcat{M}_1 := G_1 \amalg M \amalg H_1 \,,\qquad
  \Lcat{M}_0 := G_0 \amalg H_0 \,,
\end{equation*}
where coproduct $\amalg$ is the disjoint union of sets. The moment maps $l_{\Lcat{M}},r_{\Lcat{M}}: \Lcat{M}_1 \rightarrow \Lcat{M}_0$ are naturally given by the moment maps of the subsets, i.e., $l_{\Lcat{M}}(g) = l_G(g)$ and $r_{\Lcat{M}}(g) = r_G(g)$ for all $g \in G$, $l_{\Lcat{M}}(m) = l_M(m)$ and $r_{\Lcat{M}}(m) = r_M(m)$ for all $m \in M$, $l_{\Lcat{M}}(h) = l_H(h)$ and $r_{\Lcat{M}}(h) = r_H(h)$ for all $h \in H$.

\begin{Proposition}
\label{th:linking}
A set $M$ over $G_0$ and $H_0$ is a $G$-$H$ groupoid bibundle if and only if $\Lcat{M}_1 \rightrightarrows \Lcat{M}_0$ is a category which extends $G$ and $H$.
\end{Proposition}
\begin{proof}
Let $M$ be a $G$-$H$ groupoid bibundle then we can define the composition of arrows $g \in G \subset \Lcat{M}_1$ and $m \in M \subset \Lcat{M}_1$ by the left action $g \circ m := g \cdot m$ whenever $r_{\Lcat{M}}(g) = r_G(g) = l_M(m) = l_{\Lcat{M}}(m)$. Analogously, $m \circ h := m \cdot h$ whenever defined. Conversely, if $\Lcat{M}$ is a category we define the groupoid actions by the compositions, $g \cdot m := g \circ m$ and $m \cdot h := m \circ h$. Associativity of the  composition of arrows is equivalent to the the commutativity of the groupoid actions. 
\end{proof}

For a groupoid bibundle $M$ we call $\Lcat{M}$ the \textbf{linking category} associated to $M$. By definition, in $\Lcat{M}$ there are only arrows from $G_0$ to $H_0$ but no arrows from $H_0$ to $G_0$. We can express the conditions for $M$ to be right principal in terms of the linking category. Since we are working here with set-theoretic groupoids we have to view the sets as discrete, zero-dimensional manifolds. In this case, condition $(P1)$ simply states that the left moment map is surjective.

\begin{Proposition}
\label{th:linkingprincipal}
  Let $G$, $H$ be set-theoretic groupoids and $M$ a $G$-$H$ bibundle, viewed as discrete manifolds. Then the conditions for $M$ to be right principal from Definition~\ref{def:principal} can be expressed in terms of the linking category as follows:
\begin{itemize}
 \item[(P1)] $\Leftrightarrow$ For every $x \in G_0$ there is an arrow in $\Lcat{M}_1$ from $x$ to $H_0$.
 \item[(P2)] $\Leftrightarrow$ All right factorizations are unique, i.e., if for given arrows $\lambda$, $\lambda'$ there is an arrow $\mu$ such that $\lambda \circ \mu = \lambda'$ then $\mu$ is unique.
\item[(P3)] $\Leftrightarrow$ For every pair $m, m'$ of arrows from $x \in G_0$ to $H_0$ there is a right factorization $m \circ \mu = m'$.
\end{itemize}
\end{Proposition}
\begin{proof}
(P1) is obvious. 

For (P2) let $m,m' \in M \subset \Lcat{M}_1$. Then uniqueness of every right factorization $m \circ h = m \cdot h = m'$ implies that the right $H$-action is free. Conversely, if the right $H$-action is free then the right factorization of elements $m, m' \in M$ is unique. The factorizations of pairs $g,g' \in G$ and $h,h' \in H$ are unique because $G$ and $H$ are groupoids. The unique right factorization of a pair $g \in G$, $m \in M$ is given by $g \circ (g^{-1} \circ m) = m$. The opposite factorization $m \circ \mu = g$ does not exist.

For (P2) observe that every $\mu$ with $m \circ \mu = m'$ is an arrow in $H_1$. Therefore, the existence of a factorization is equivalent to $H$ acting transitively on the left fibers of $M$.
\end{proof}

If $M$ is biprincipal we can go further and define a groupoid associated to $M$. First, note that for every $G$-$H$ bibundle $M$ there is an opposite $H$-$G$ bibundle $M^\op$ given by the set $M$ with swapped moment maps $l_{M^\op}(m) = r_M(m)$, $r_{M^\op}(m) = l_M(m)$ and actions $h \cdot_{\op} m = m \cdot h^{-1}$, $m \cdot_{\op} g = g^{-1} \cdot m$. 

\begin{Proposition}
\label{th:linkinggroupoid}
Let $G$, $H$ be set-theoretic groupoids and $M$ a biprincipal $G$-$H$ bibundle. Consider the sets
\begin{equation*}
  \Lgpd{M}_1 := G_1 \amalg M \amalg M^\op \amalg H_1 \,,\qquad
  \Lgpd{M}_0 := G_0 \amalg H_0
\end{equation*}
with left and right moment maps given by the left and right moment maps of each summand and with the compositions
\begin{equation*}
\begin{gathered}
  g \circ g' := gg' \,,~
  m \circ \bar{m}' := {}_G\langle m, \bar{m}' \rangle \,,
\\
  \bar{m} \circ g' := \bar{m} \cdot_\op g' \,,~
  h \circ \bar{m}' := h \cdot_\op \bar{m}' \,,
\end{gathered}
\quad
\begin{gathered}
  g \circ m' := g \cdot m' \,,~
  m \circ h := m \cdot h' \,,
\\
  \bar{m} \circ m' := \langle \bar{m}, m' \rangle_H \,,~
  h \circ h' := hh' \,,
\end{gathered}
\end{equation*}
whenever defined for $g,g \in G_1$, $m,m' \in M$, $\bar{m}, \bar{m}' \in M^\op$, $h,h' \in H$, where ${}_G\langle ~, ~ \rangle$ and $\langle ~, ~ \rangle_H$ denote the $G$-valued and the $H$-valued bibundle pairings, respectively.

Then $\Lgpd{M}_1 \rightrightarrows \Lgpd{M}_0$ is a groupoid, called the \textbf{linking groupoid} of $M$.
\end{Proposition}
\begin{proof}
Associativity of the composition follows from the commutativity of the groupoid actions and the properties of the bibundle pairings. The inverse of $m \in M$ is $m \in M^\op$.
\end{proof}

If $G$ and $H$ are Lie groupoids and $M$ a smooth biprincipal bibundle then $\Lgpd{M}$ is a Lie groupoid. The linking groupoid is then the differential geometric analog of the linking algebra of a Hilbert bimodule \cite{BrownGreenRieffel:Stable}. In fact, the order in which the groupoid composition was written above is to resemble the ``multiplication'' of matrices $(\begin{smallmatrix} g & m \\ \bar{m} & h \end{smallmatrix})$ and $(\begin{smallmatrix} g' & m' \\ \bar{m}' & h' \end{smallmatrix})$ with the sums omitted. Apparently, the concept of linking groupoids has not appeared in the literature but was known to some experts.\footnote{Andr\'e~Haefliger, private communication.}

\subsection{Principality as a categorical property}

In the category of Lie groupoids and smooth homomorphisms of groupoids finite products and a terminal object exist. For two groupoids $G$ and $H$ the categorical product is given by the cartesian product $G \times H$ together with the cartesian product of all the structure maps. The terminal object $1$ is the groupoid with one object and one arrow. The diagonal map $G \rightarrow G \times G$, $g \mapsto (g,g)$ and the terminal map $G \rightarrow 1$, $g \mapsto 1$ are homomorphisms of groupoids.

\begin{Definition}
For groupoids $G$ and $H$ denote the bundlization of the diagonal map $G \rightarrow G\times G$ by $\Delta_G$, and the bundlization of the terminal map $G \rightarrow 1$ by $\varepsilon_G$.
\end{Definition}

Explicitly, $\Delta_G$ is the set of all arrows with the same left end,
\begin{equation*}
  \Delta_G := \{ (g_1,g_2) \in G_1 \,|\, l(g_1) = l(g_2) \} \,,
\end{equation*}
with the left action $g' \cdot (g_1,g_2) = (g'g_1,g'g_2)$ and the right action $(g_1,g_2)\cdot (g'_1,g'_2) = (g_1 g'_1,g_2 g'_2)$ whenever defined for $(g_1,g_2) \in \Delta_G$, $g' \in G$ and $(g'_1,g'_2) \in G \times G$. The bibundle $\varepsilon_G$ is the base $G_0$ of the groupoid with the left action of $G$ on its base, $g \cdot x = l(g)$.

The following theorem shows that right principality of a bibundle is a categorical property: 
\begin{Theorem}
\label{th:rightprincipal}
Let $G$ and $H$ be Lie groupoids and $M$ be a smooth $G$-$H$-bibundle. $M$ is right principal if and only if
\begin{equation}
\label{eq:principaldiag}
\begin{gathered}
\xymatrix@C-2ex{
  G \ar[rr]^{M} \ar[dr]_{\varepsilon_G} 
  & \dtwocell<\omit>{} & H \ar[dl]^{\varepsilon_H} \\
  & 1 }
\qquad\qquad
\xymatrix{
G \ar[r]^{M} \ar[d]_{\Delta_G} \drtwocell<\omit>{} & H \ar[d]^{\Delta_H}\\
G \times G \ar[r]^{M \times M} & H \times H }
\end{gathered}
\end{equation}
are 2-commutative diagrams of smooth bibundles.
\end{Theorem}

We will prove it in two steps, by showing what the weak commutativity of each one of the two diagrams means geometrically.

\begin{Proposition}
\label{th:righttransitive}
Let $G$, $H$ be set-theoretic groupoids and $M$ a $G$-$H$-bibundle. $M$ satisfies properties (P1) and (P3) of Definition~\ref{def:principal} if and only if the left diagram of~\eqref{eq:principaldiag} is 2-commutative.
\end{Proposition}

\begin{Lemma}
\label{th:terminal1}
  Let $M$ be a $G$-$1$ bibundle. The left moment map $l_M : M \rightarrow G_0 = \varepsilon_G$ is the unique $G$-equivariant map from $M$ to $\varepsilon_G$. If $M$ is right principal, then $l_M$ is an isomorphism of bibundles.
\end{Lemma}

\begin{proof}
Let $M$ be a $G$-$1$ bibundle. The left moment map $l_{M} : M \rightarrow G_0 = \varepsilon_G$ satisfies $l_{M}(g \cdot m) = l_G(g) = g \cdot l_{M}(m)$, so it is equivariant. Assume that $\phi: M \rightarrow \varepsilon_G$ is another equivariant map. Because $l_{\varepsilon_G} = \id_{G_0}$ the equivariance of $\phi$ implies $l_{M} = l_{\varepsilon_G} \circ \phi = \id_{G_0} \circ \phi = \phi$. We conclude that $l_{M}$ is the only equivariant map from $l_{M}$ to $\varepsilon_{G}$.

Note, that $\varepsilon_{G}$ is trivially right principal. Let now $M$ be another right principal $G$-$1$ bundle. Since the terminal groupoid $1$ acts trivially, every left fiber $l_{M}^{-1}(x)$, $x \in G_0$ consists of a single point. Therefore, the equivariant map $l_{M}: M \rightarrow \varepsilon_{G}$ is invertible, the inverse being automatically equivariant.
\end{proof}

\begin{proof}[Proof of Prop.~\ref{th:righttransitive}]
Consider the $G$-$1$ bibundle $N := M \circ \varepsilon_H$. By the last lemma the left moment map $l_N: N \rightarrow \varepsilon_G$ is the unique equivariant map. As set $N = M/H$, the set of orbits of the right $H$-action on $M$. The left moment map $l_N$ maps $[m] \in M/H$ to $l_M(m)$. This map is surjective iff $l_M$ is surjective. It is injective iff there is only one orbit in each fiber, which is the case iff the right action of $H$ is transitive on the left fibers.
\end{proof}

This proposition gives a categorical definition for a bibundle $M$ to be transitive with surjective left moment map and a categorical description of the right orbit space as composition of bibundles $M \circ \varepsilon_H$.

\begin{Corollary}
If a $G$-$H$ bibundle $M$ and a $H$-$K$ bibundle $N$ are both right transitive (on left fibers), then $M \circ N$ is also right transitive.
\end{Corollary}

\begin{Proposition}
\label{th:principal2}
Let $G$, $H$ be set-theoretic groupoids and $M$ a $G$-$H$-bibundle. $M$ satisfies properties (P2) and (P3) of Definition~\ref{def:principal} if and only if the right diagram of~\eqref{eq:principaldiag} is 2-commutative.
\end{Proposition}

\begin{proof}
According to Proposition~\ref{th:freetransitive} $M$ satisfies properties (P2) and (P3) if and only if there is a bibundle pairing. Given a pairing, consider the maps
\begin{xalignat*}{2}
  M \times_{G_0} M
  &\longrightarrow M \times_{H_0} H
  &
  M \times_{H_0} H
   &\longrightarrow M \times_{G_0} M 
\\
  (m,m') &\longmapsto (m, \langle m , m' \rangle)
  &
  (m,h) &\longmapsto (m , m \cdot h ) 
\end{xalignat*}
which are mutually inverse. Both maps become $G$-$(H \times H)$ biequivariant, if we equip $M \times_{G_0} M$ with the left moment map $l(m,m') := l(m) = l(m')$, with the right moment map $r(m,m') = (r(m),r(m'))$, and with the actions
\begin{equation*}
  g \cdot (m,m') = (g \cdot m, g \cdot m') \,,\qquad
  (m,m') \cdot (h,h') = (m \cdot h, m \cdot h') \,,
\end{equation*}
and if we equip $M \times_{H_0} H$ with the left moment map $l(m,h) = l(m)$, with the right moment map $r(m,g) = (r(m),r(g))$, and with the left and right actions
\begin{equation*}
  g \cdot (m,h) = (g \cdot m, h) \,,\qquad
  (m,h) \cdot (h_1,h_2) = (m \cdot h_1, h_1^{-1} g h_2) \,.
\end{equation*}
There is an isomorphism of bibundles 
\begin{equation*}
  M \circ \Delta_H
  \stackrel{\cong}{\longrightarrow}
  M \times_{H_0}^{r,l} H 
  \,,\qquad
  [m,h,h'] \longmapsto (m\cdot h, h^{-1} h') \,,
\end{equation*}
the inverse map being given by $(m,h) \mapsto [m,1_{r(m)},h]$. And there is an isomorphism of bibundles
\begin{equation*}
  \Delta_G \circ (M \times M)
  \stackrel{\cong}{\longrightarrow}
  M \times_{G_0}^{l,l} M 
  \,,\qquad
  [g,m,m'] \longmapsto (g \cdot m, g \cdot m') \,,
\end{equation*}
the inverse being given by $(m,m') \mapsto [1_{l(m)},m,m']$. We conclude that (P2) and (P3) are satisfied if and only if there is an isomorphism of bibundles $M \circ \Delta_H \cong \Delta_G \circ (M \times M)$.
\end{proof}

\begin{proof}[Proof of Theorem~\ref{th:rightprincipal}]
Let $M$ be a smooth bibundle such that the diagrams~\eqref{eq:principaldiag} smoothly 2-commute. Prop.~\ref{th:principal2} implies that $M$ satisfies (P2) and (P3). Prop.~\ref{th:righttransitive} implies that the left moment map of $M$ is surjective. The left moment map of $\varepsilon_G \cong M \circ \varepsilon_H$ is submersive, which implies that the left moment map of $M \times_{H_0} \varepsilon_H$ and, thus, the left moment map of $M$ is submersive. We conclude that $M$ is smooth right principal. 

Conversely, if $M$ is smooth right principal, then all compositions of bibundles in the diagrams are smooth and Prop.~\ref{th:righttransitive} and Prop.~\ref{th:principal2} imply that the diagrams are 2-commutative. The 2-isomorphism which makes the first diagram 2-commutative is given by the left moment map of $M \circ \varepsilon_H$ which is smooth because $M$ is smooth. The 2-isomorphism which makes the second diagram two-commutative is essentially given by the bibundle pairing which was shown to be smooth in Proposition~\ref{th:freetransitive}. We conclude that both diagrams~\eqref{eq:principaldiag} are smooth 2-commutative.
\end{proof}

\begin{Corollary}
\label{th:principalcomp}
If a smooth $G$-$H$ bibundle $M$ and a smooth $H$-$K$ bibundle $N$ are both right principal, then $M \circ N$ is also smooth right principal.
\end{Corollary}

The bibundles $\Delta_G$ and $\varepsilon_G$ are defined as bundlizations of the diagonal and the terminal groupoid homomorphism of $G$, respectively. The bundlization of the identity map of $G$ is the groupoid itself, $\Id_G = G$. Furthermore, we have natural 1-isomorphisms $1 \times G \cong G \cong G \times 1$. By functoriality of bundlization we obtain that the diagrams
\begin{equation}
\label{eq:AssocUnitalDiag}
\begin{gathered}
\xymatrix@+1ex{
G  
  \ar[r]^-{\Delta \times \Id} \ar[d]_{\Id \times \Delta_G}
  \drtwocell<\omit>{}
& G\times G 
  \ar[d]^{\Delta_G} 
\\ 
G \times G \ar[r]^{\Delta_G} & G \times G \times G
}
\qquad
\xymatrix@+1ex{
 & G \ar[dl]_{\cong} 
     \ar[d]^{\Delta_G} 
     \ar[dr]^{\cong} \drtwocell<\omit>{<2>}\\
1 \times G  \urtwocell<\omit>{<2>} & G \times G 
  \ar[l]^-{\varepsilon_G \times \Id_G} 
  \ar[r]_-{\Id_G \times \varepsilon_G} & G \times 1
}
\end{gathered}
\end{equation} 
are 2-commutative diagrams of smooth bibundles. In this sense, the diagonal and terminal bibundle equip the groupoid with the structure of a weak comonoid. Theorem~\ref{th:rightprincipal} can than be rephrased by saying that a bibundle is right principal if and only if it is a weakly comonoidal 1-morphism.

\subsection{Products}

In the weak 2-category of set-theoretic groupoids and bibundles (not necessarily right principal) every left $G$-bundle is a morphism from $G$ to the trivial groupoid $1$. We conclude that in this category $1$ is not a terminal object. But it has been shown in Lemma~\ref{th:terminal1} that if $T$ is a \emph{right principal} $G$-$1$ bibundle, then the left moment map $l_T$ is the unique isomorphism of bibundles from $T$ to $\varepsilon_G$. In this sense, $1$ is a terminal object in the weak 2-category of groupoids and right principal bibundles. Rigorously speaking, $1$ is a terminal 2-categorical cone over the empty functor. Similarly, the product groupoid $G \times H$ is not the 2-categorical product in the category of groupoids and bibundles. We will now show that right principality is a sufficient and necessary condition which the bibundles in a weak 2-subcategory have to satisfy for $G \times H$ to be the 2-categorical product. For limits in higher categories we refer the reader to \cite{Borceux:Handbook1} \cite{Joyal:Quasi} and section~\ref{sec:GroupsHigher}.

\begin{Proposition}
\label{th:CartesianCat}
  Let $G$, $H$, $K$ be Lie groupoids. Let $\pr_G := \Id_G \times \varepsilon_H$ and $\pr_H := \varepsilon_G \times \Id_H$. For every smooth right principal $K$-$G$ bibundle $M$ and every smooth right principal $K$-$H$ bibundle $N$, the smooth right principal $K$-$(G\times H)$ bibundle $L := \Delta_K \circ (M \times N)$ makes the diagram
\begin{equation}
\label{eq:ProductDiag}
\begin{gathered}
  \xymatrix@+1ex{
 & K \ar[dl]_{M} \ar[d]^{L} \ar[dr]^{N} \drtwocell<\omit>{<2>}\\
G  \urtwocell<\omit>{<2>} & G\times H \ar[l]^-{\pr_G} \ar[r]_-{\pr_H} & H
}
\end{gathered}
\end{equation}
2-com\-mu\-ta\-tive. Every other right principal bibundle which makes the diagram 2-com\-mu\-ta\-tive is isomorphic to $L$.
\end{Proposition}

\begin{proof}
By Proposition~\ref{th:righttransitive} and the commutativity of the second of diagrams~\eqref{eq:AssocUnitalDiag} we get
\begin{equation*}
\begin{split}
  L \circ \pr_G 
  &\cong \Delta_K \circ (M \times N)
     \circ(\Id_G \times \varepsilon_H)
  \cong \Delta_K \circ (M \times [N \circ \varepsilon_H]) \\
  &\cong \Delta_K \circ (M \times \varepsilon_K) 
  \cong M \,,
\end{split}
\end{equation*}
and analogously for the projection onto $H$. We conclude that diagram~\eqref{eq:ProductDiag} is weakly commutative.

It remains to show that $L$ is unique up to 2-isomorphism. Let $L'$ be another right principal $K$-$(G\times H)$ bibundle which makes the diagram commutative. Because $L'$ is right principal we get from Theorem~\ref{th:rightprincipal} an isomorphism of bibundles,
\begin{equation*}
  L' \circ \Delta_{G \times H}
  \stackrel{\cong}{\longrightarrow} 
  \Delta_K \circ (L' \times L') \,. 
\end{equation*}
We can compose this from the right with the identity map on the bibundle $\pr_G \times \pr_H$ and obtain an isomorphism of bibundles
\begin{equation}
\label{eq:Product1}
  L' \circ \Delta_{G \times H} 
  \circ (\pr_G \times \pr_H)
  \stackrel{\cong}{\longrightarrow} 
  \Delta_K \circ (L' \times L') 
  \circ (\pr_G \times \pr_H) \,. 
\end{equation}
From the definition of the projection bibundles and the commutativity of diagrams~\eqref{eq:AssocUnitalDiag} we obtain for the left hand side
\begin{equation*}
  L' \circ \Delta_{G \times H} \circ (\pr_G \times \pr_H) 
  \cong L' \,.
\end{equation*}
The assumption that $L'$ makes the diagram~\eqref{eq:ProductDiag} weakly commute implies for the right hand side of Eq.~\eqref{eq:Product1} that
\begin{equation*}
  \Delta_K \circ (L' \times L') 
  \circ (\pr_G \times \pr_H) \cong
  \Delta_K \circ (M \times N) 
  \cong L \,.
\end{equation*} 
From the last three equations we we conclude that $L' \cong L$.
\end{proof}

This shows that right principality of bibundles is sufficient for $G \times H$ to be the weak categorical product in the sense of the last proposition. It is also a necessary condition:

\begin{Proposition}
\label{th:CartesianCat2}
 Let $\calC$ be a weak 2-subcategory of the category of set-theoretic group\-oids and groupoid bibundles (not necessarily right principal). If for all groupoids $G$ and $H$ the cartesian product $G \times H$ is the weak categorical product in the sense of Proposition~\ref{th:CartesianCat} then all bibundles in $\calC$ are right principal.
\end{Proposition}

\begin{proof}
The universal property of the categorical product implies that there are bibundles which are isomorphic to $\varepsilon_G$ and $\Delta_G$. By abuse of notation, we denote these with the same letters. Let $M$ be a $K$-$G$ bibundle. By the universal property, the diagrams
\begin{equation*}
\begin{gathered}
  \xymatrix@+1ex{
 & K \ar[dl]_{M} \ar[d]^{L} \ar[dr]^{\varepsilon_K} \drtwocell<\omit>{<2>}\\
1  \urtwocell<\omit>{<2>} & G \times 1 \ar[l]^-{\pr_G} \ar[r]_-{\pr_1} & 1
}
\qquad
  \xymatrix@+1ex{
 & K \ar[dl]_{M} \ar[d]^{L} \ar[dr]^{M} \drtwocell<\omit>{<2>}\\
1  \urtwocell<\omit>{<2>} & G \times G \ar[l]^-{\pr_G} \ar[r]_-{\pr_G} & G
}
\end{gathered}
\end{equation*}
are 2-commutative. The left diagram tells us that $L \circ \pr_G \cong M$. Because $G \times 1 \cong G$ we see that $\pr_G \cong \Id_G$ so $L \cong M$. Similarly, $\pr_1 \cong \varepsilon_G$. With this the 2-commutativity of the left triangle then implies that $M \circ \varepsilon_G \cong \varepsilon_K$. The right diagram implies that $L :=  \Delta_K \circ (M \times M)\cong M \circ \Delta_G$. From Theorem~\ref{th:rightprincipal} it follows that $M$ is right principal. 
\end{proof}

\subsection{The groupoid inverse as rigid structure}

As for groups, the inverse of a groupoid can be viewed as a natural transformation from $G$ to $G^\op$, where $G^\op$ is the groupoid with the same sets of arrows and the same base, $G_1^\op := G_1$, $G_0^\op := G_0$, swapped moment maps $l^\op := r$, $r^\op = l$, and opposite multiplication $g \cdot_\op g' = g'g$. The inverse can then be used to convert a left into a right groupoid action, and vice versa. To every smooth $G$-$H$ bibundle $M$ we can thus naturally associate a $G^\op$-$H^\op$ bibundle, which is the same as a $H$-$G$ bibundle. We will call this the opposite bibundle $M^\op$. In other words, we have a natural bijection from $\Hom_1(G,H)$ to $\Hom_1(H,G)$ in $\GrpdBibu$. 

But we can also convert only one of the actions of a bibundle. In this way we obtain from $M$ a $G$-$H^\op$ bibundle, which we can view as $(G \times H)$-$1$ bibundle or as $(H \times G)$-$1$ bibundle. We can avoid this ambiguity by describing the conversion of right into left actions also as the composition of bibundles:  Consider the $(G \times G)$-$1$ bibundle $\Ev_G$, called evaluation, and the $1$-$(G \times G)$ bibundle $\Cv_G$, called coevaluation, defined as
\begin{xalignat}{3}
\label{eq:rigidinverse2a}
  \Ev_G &= G\,, & l_{\Ev}(g) &= (l_G(g),r_G(g)) \,, 
    & (g_1,g_2) \cdot g &= g_1 g g_2^{-1} \,,\\
\label{eq:rigidinverse2b}
  \Cv_G &= G\,, & r_{\Cv}(g) &= (r_G(g),l_G(g)) \,, 
    & g \cdot (g_1,g_2) &= g_1^{-1} g g_2 \,.
\end{xalignat}
Since the groupoid inverse is an involutive operation, the diagram
\begin{equation}
\label{eq:RigidDiag}
\begin{gathered}
\xymatrix@d@R-7ex@C-2ex{
  1 \times G \ar[r]^{\Cv_G \times \Id_G}
  \xtwocell[rrr]{}<\omit>{<6>}
& (G \times G) \times G \ar[r]^{\cong}
& G \times (G \times G) \ar[r]^{\Id_G \times \Ev_G}
& G \times 1 \ar[d]^-{\cong}
\\
  G \ar[u]^-{\cong} \ar[d]_-{\cong}
  \ar[rrr]^{\Id_G}
&&& G
\\
  G \times 1 \ar[r]_{\Id_G \times \Cv_G}
& G \times (G \times G) \ar[r]_{\cong}
& (G \times G) \times G \ar[r]_{\Ev_G \times \Id_G}
& 1 \times G \ar[u]_-{\cong}
  \xtwocell[lll]{}<\omit>{<6>}
}
\end{gathered} 
\end{equation} 
is a 2-commutative diagram of 1-morphisms in $\GrpdBibu$. Up to 2-isomorphisms this is the diagram of a rigid structure of the monoidal category given by the cartesian product of groupoids and the trivial groupoid, where we have to define the (right and left) rigid dual simply by 
\begin{equation*}
  G^\vee := G \,.
\end{equation*}
Let us introduce the notation $\GrpdBibu(1)$ to denote the category associated to the bicategory $\GrpdBibu$ which we obtain by modding out 2-isomorphism. Now we can summarize this observation rigorously as follows:

\begin{Proposition}
\label{th:rigid}
The category $\GrpdBibu(1)$ together with the cartesian product as tensor product, the trivial groupoid as unit object, and the coevaluation and evaluation 1-morphisms defined as in Eqs.~\eqref{eq:rigidinverse2a} and \eqref{eq:rigidinverse2b} is a rigid monoidal category.
\end{Proposition}

Using the rigid structure we can express the opposite of a $G$-$H$ bibundle up to 2-isomorphism as
\begin{equation}
\label{eq:BibuOpRigid}
  M^{\op} \cong
  \bigl( (\Cv_G \times \Id_H) \circ 
  (\Id_G \times M \times \Id_H) \bigr) \circ
  (\Id_G \times \Ev_H ) \,.
\end{equation}
Note, that the coevaluation and evaluation bibundles are not right principal, so this rigid structure does not live in $\LieGrpdPrBibu$. If we compose the rigid structure for smooth bibundles we have to decide on a case to case basis whether the resulting bibundle is smooth, or whether we have to work with topological or set-theoretic groupoids. The coevaluation bibundle satisfies the following useful property: 

\begin{Proposition}
The diagram
\begin{equation}
\label{eq:CvDelta}
\begin{gathered}
\xymatrix@R-2ex@C-3ex{
 & 1
  \ar[dl]_-{\Cv_G}
  \ar[dr]^-{\Cv_G}
  \ar[ddd]^{\Cv_G}
 \xtwocell[ddd]{}<\omit>{<-6>} 
\\
G \times G
  \ar[d]_{\Delta_G \times \varepsilon_G } 
&& G \times G
  \ar[d]^{\varepsilon_G \times \Delta_G} 
\\
(G \times G) \times 1
  \ar[dr]_{\cong}
&& 1 \times (G \times G)
  \ar[dl]^{\cong}
\\
& G \times G
 \xtwocell[uuu]{}<\omit>{<-6>} 
}
\end{gathered}
\end{equation}
is a 2-commutative diagram of smooth bibundles.
\end{Proposition}

\subsection{The graphic notation}

As mentioned in the last section, the coherence axioms of a weak 2-category ensure that we obtain an honest category when we mod out the 2-morphisms and consider isomorphism classes of bibundles rather than the bibundles themselves.  Isomorphism classes of right principal bibundles are often called Hilsum-Skandalis morphisms and have first appeared in the theory of foliations \cite{HilsumSkandalis:Morphismes}.

Consider the full monoidal subcategory generated by a groupoid $G$. The objects in this subcategory are the products of $G$, $G^0 = 1$, $G^1 = G$, $G^2 = G \times G$, etc. The morphisms are $G^m$-$G^n$-bibundles for $m,n \in \bbN_0$. By definition of the composition of bibundles and because we are in a symmetric monoidal category we can selectively compose every single factor of $G$ acting on the right on a bibundle $M$ in this subcategory with every factor of $G$ acting on the left on another bibundle $N$. On isomorphism classes of bibundles this can be extended to the structure of a PROP. Here, we will not make use of the full structure of a PROP but only use the usual graphic notation:

A morphism given by the isomorphism class $[M]$ of a $G^m$-$G^n$-bibundle is represented by a box with $m$-lines entering on the top and $n$-lines exiting at the bottom. (The usage of square brackets has nothing to do with the functor from groupoids to stacks in Sec.~\ref{sec:GrpdsToStacks}.) The cartesian product of two such morphisms is represented by horizontal juxtaposition.  The composition of morphisms is represented by stacking the morphisms from top to bottom, connecting the lines in the middle. For example, if $[M]$ is the isomorphism class of a $G^2$-$G$-bibundle and $[N]$ of a $G$-$G^2$-bibundle we write
\begin{equation*}
[M] \times [N] \quad = \quad
\begin{xyc}
  \MorpM(-1,0) \MorpN(2,0)
\end{xyc}
\qquad\text{and}\qquad
[M] \circ [N] \quad = \quad
\begin{xyc}
  \MorpM(0,0) \MorpN(0,-2)
\end{xyc} \,.
\end{equation*}
Because we are now working with isomorphism classes the composition is strictly associative.

We have already obtained a number of bibundles from the basic structure of the groupoid categories: the terminal object, the identity, the bibundle associated to the diagonal homomorphism and the terminal homomorphism, the evaluation and coevaluation bibundles. The categorical product is symmetric so we have a flip $\tau_{GH}$, given by the isomorphism class of the bundlization of $\tau: G \times H \rightarrow H \times G$, $(g,h) \mapsto (h,g)$. Because these bibundles play a fundamental role, we simplify their notation and suppress the box and the letter denoting the bibundle: 
\begin{alignat*}{5}
&\text{terminal object:}& 1 &=
\begin{xyc}
  \Term(0,0)
\end{xyc}
&&\text{flip:}&& [\tau_{GG}] &=
\begin{xyc}
  \Flip(0,0)
\end{xyc}\\[1ex]
&\text{identity of $G$:}& [\Id_G] &=
\begin{xyc}
  \Iden(0,0)\Iden(0,-1)
\end{xyc}
&&\text{evaluation:}&& [\Ev_G] &=
\begin{xyc}
  \Coev(0,0)
\end{xyc}
\\[1ex]
&\text{diagonal:} & [\Delta_G] &=
\begin{xyc}
  \Prod(0,0)
\end{xyc}
\qquad\quad
&& \text{coevaluation:~}&& [\Cv_G] &= 
\begin{xyc}
  \Eval(0,0)
\end{xyc}
\\[1ex]
&\text{terminal morphism:} & [\varepsilon_G] &=
\begin{xyc}
  \Unit(0,0)
\end{xyc}
& && && & \quad.
\end{alignat*}
Because we are working in the monoidal subcategory generated by $G$ we will also drop the subscript denoting the groupoid, $\Id_G = \Id$, etc. 

When expressed in this graphic notation, the properties of these bibundles take on a very intuitive form. For example, the property of the identity simply means that every vertical line can be extended or shortened. The weak commutativity of diagrams~\eqref{eq:AssocUnitalDiag} states that the diagonal and the terminal bibundle are a comonoidal structure:
\begin{equation}
\label{eq:comonoidal}
\begin{xyc}
  \Prod(0,0) \Prod(1,2) \Iden(2,-1)\Iden(2,0)
\end{xyc}
\quad = \quad
\begin{xyc}
  \Prod(0,0) \Prod(-1,2) \Iden(-2,-1)\Iden(-2,0)
\end{xyc} 
\qquad\text{and}\qquad
\begin{xyc}
  \Prod(0,1) \Unit(-1,-1) \Iden(1,-1)
\end{xyc}
\quad = \quad
\begin{xyc}
  \Iden(0,-1)\Iden(0,0)\Iden(0,1)
\end{xyc}
\quad = \quad
\begin{xyc}
  \Prod(0,1) \Unit(1,-1) \Iden(-1,-1)
\end{xyc} \quad.
\end{equation}
Theorem~\ref{th:rightprincipal} states that a $G$-$G$-bibundle $M$ is right principal if and only if: 
\begin{align*}
\begin{xyc}
  \Unit(0,0) \Morp(0,1){M}
\end{xyc}
\quad = \quad
\begin{xyc}
  \Unit(0,0) \Iden(0,1)
\end{xyc}
\qquad\text{and}\qquad
\begin{xyc}
  \Morp(0,0){M} \Prod(0,-1)
\end{xyc}
\quad = \quad
\begin{xyc}
  \Prod(0,0) \Morp(-1,-3){M} \Morp(1,-3){M}
\end{xyc}
 \quad.
\end{align*}
The weak commutativity of diagram~\eqref{eq:RigidDiag} states that $[\Ev_G]$ and $[\Cv_G]$ are a rigid structure:
\begin{equation*}
\begin{xyc}
  \Coev(0,0) \Iden(3,0) \Iden(-1,1) \Eval(2,2)
\end{xyc}
\quad = \quad
\begin{xyc}
  \Iden(0,0)\Iden(0,1)
\end{xyc}
\quad = \quad
\begin{xyc}
  \Coev(2,0) \Iden(3,1) \Iden(-1,0) \Eval(0,2)
\end{xyc}
\end{equation*}
The weak commutativity of diagram~\eqref{eq:CvDelta} is expressed as:
\begin{equation}
\label{eq:Gcoevdiag}
\begin{xyc}
  \Prod(0,0)\Eval(1,2) \Unit(2,0)
\end{xyc}
\quad = \quad
\begin{xyc}
  \Eval(0,1) \Iden(-1,-1) \Iden(1,-1)
\end{xyc}
\quad = \quad
\begin{xyc}
  \Prod(0,0)\Eval(-1,2) \Unit(-2,0)
\end{xyc}
\end{equation}
Eq.~\eqref{eq:BibuOpRigid} shows how the opposite of a $G$-$G$-bibundle $M$ is expressed in terms of the rigid structure:
\begin{equation*}
[M^{\op}]
\quad = \quad
\begin{xyc}
  \Eval(-1,3) \Morp(0,0){M} \Coev(1,-1)
  \Iden(-2,1)\Iden(-2,0)
  \Iden(2,0)\Iden(2,1)
\end{xyc} \,.
\end{equation*}
The left version of Prop.~\ref{th:righttransitive} states that $M$ is left transitive and right surjective if and only if:
\begin{equation*}
\begin{xyc}
  \Eval(1,3) 
  \Morp(0,0){M} \Unit(2,1)
\end{xyc}
\quad\cong\quad 
\begin{xyc}
  \Eval(1,3) 
  \Iden(0,1) \Unit(2,1)
\end{xyc}
\quad.
\end{equation*}
The flip is idempotent, that is, the monoidal structure is symmetric, and natural,
\begin{equation*}
\begin{xyc}
  \Flip(0,0)\Flip(0,2)
\end{xyc}
\quad = \quad
\begin{xyc}
  \Iden(-1,-1)\Iden(-1,0)\Iden(-1,1)\Iden(-1,2)
  \Iden(1,-1)\Iden(1,0)\Iden(1,1)\Iden(1,2)
\end{xyc}
\qquad\text{and}\qquad
\begin{xyc}
  \Flip(0,0)\Morp(-1,1){M}\Morp(1,1){N}
\end{xyc}
\quad = \quad
\begin{xyc}
  \Flip(0,2)\Morp(-1,-1){N}\Morp(1,-1){M}
\end{xyc}
\quad,
\end{equation*}
for all $G$-$G$-bibundles $M$,$N$. The diagonal bibundle is symmetric:
\begin{equation*}
\begin{xyc}
  \Flip(0,2)\Prod(0,4)
\end{xyc}
\quad = \quad
\begin{xyc}
  \Iden(-1,2)\Iden(1,2)\Prod(0,4)
\end{xyc}
\quad.
\end{equation*}
From this and diagram~\eqref{eq:Gcoevdiag} it follows that the evaluation and coevaluation bibundles are symmetric, too:
\begin{equation*}
\begin{xyc}
  \Flip(0,2)\Eval(0,4)
\end{xyc}
\quad = \quad
\begin{xyc}
  \Iden(-1,2)\Iden(1,2)\Eval(0,4)
\end{xyc}
\qquad\text{and}\qquad
\begin{xyc}
  \Flip(0,2)\Coev(0,0)
\end{xyc}
\quad = \quad
\begin{xyc}
  \Iden(-1,2)\Iden(1,2)\Coev(0,1)
\end{xyc} \quad.
\end{equation*}
The graphic notation invites to tinker with the given bibundles. For example, here are two bibundles and their geometric meaning:
\begin{align*}
\begin{xyc}
  \Eval(0,2) \Unit(-1,0) \Unit(1,0)
\end{xyc}
\quad &= \quad \text{[coarse moduli space of $G$]}\quad = \quad G_0/G_1 \,,\\[1ex]
\begin{xyc}
  \Prod(0,0) \Coev(0,-2)
\end{xyc}
\quad &= \quad \text{[isotropy bundle of $G$ with left adjoint action.]}
\end{align*}
This representation of the coarse moduli space is an example for the fact that compositions with the coevaluation bibundle, which is not right principal, may not be smooth. As illustration for the usefulness of the graphic notation we give a concise graphic proof of the following well known fact:

\begin{Proposition}
  Let $G$ be a Lie groupoid and $M$ a smooth $G$-$G$ bibundle. If $M$ is biprincipal then $M^{\op}$ is a weak 1-inverse.
\end{Proposition}
\begin{proof}
Assume that $M$ is biprincipal. Then
\begin{equation*}
\begin{split}
  [M^{\op}] \circ [M] 
  &=\quad
\begin{xyc}
  \Eval(-1,3) \Morp(0,0){M} \Coev(1,-1)
  \Morp(-2,0){M}
  \Iden(2,0)
\end{xyc}
\quad = \quad
\begin{xyc}
  \Eval(0,5)
  \Unit(1,3)
  \Prod(-1,3) \Morp(0,0){M} \Coev(1,-1)
  \Morp(-2,0){M}
  \Iden(2,0)
\end{xyc}
\quad = \quad
\begin{xyc}
  \Eval(0,5)
  \Unit(1,3)
  \Morp(-1,2){M}
  \Prod(-1,1) \Coev(1,-1)
  \Iden(2,0)
\end{xyc}
\\
&= \quad
\begin{xyc}
  \Eval(0,5)
  \Unit(1,3)
  \Prod(-1,3) \Coev(1,1)
  \Iden(-2,1)\Iden(2,2)
\end{xyc}
\quad = \quad
\begin{xyc}
  \Eval(-1,3) \Coev(1,1)
  \Iden(-2,1)\Iden(2,2)
\end{xyc}
\quad = \quad
\begin{xyc}
  \Iden(0,1)\Iden(0,2)
\end{xyc}
\quad.
\end{split}
\end{equation*}
Hence, $[M^\op]$ is the left inverse of $[M]$. The proof that $[M^\op]$ is the right inverse is analogous.
\end{proof}

This proof, and in fact the whole graphic notation, can be generalized to a $G$-$H$ bibundle by coloring the lines in order to distinguish the $G$ slots from the $H$ slots. We will also need the converse statement, for which we do not know a graphic proof:

\begin{Proposition}
\label{th:isoprincipal1}
  Let $G$ be a Lie groupoid, and $M$ a smooth $G$-$H$ bibundle. If $M$ is a weak 1-isomorphism (a Morita equivalence) then $M$ is biprincipal.
\end{Proposition}
\begin{proof}
Assume that given a smooth $G$-$H$ bibundle $M$ there is a smooth $H$-$G$ bibundle such that $M \circ N \cong \Id_G$ and $N \circ M \cong \Id_H$. Then $M \circ N$ is right principal, so the moment map $l_{M \circ N}$ which maps $[m,n] \in M \circ N$ to $l_{M}(m)$ is surjective, which implies that $l_{M}$ must be surjective. By analogous arguments we can show that $r_{M}$, $l_{N}$, and $r_{N}$ are surjective.

Assume the right action of $G$ on $N$ were not free. Then there is a $g \in G$ such that $n \cdot g = n$. Because $r_{M}$ is surjective there is an $m \in M$ with $r_{M}(m) = l_{N}(n)$, so there is an element in $M \circ N$ represented as $[m,n]$. Then $[m,n] \cdot g = [m,n \cdot g] = [m,n]$, so the right action of $G$ on $M \circ N$ cannot be free either, which is a contradiction. By analogous arguments we can show that $N$ and $M$ are left and right free.

Let $n,n' \in N$ be elements of the same left fiber, $l_{N}(n) = l_{N}(n')$.  Because $r_M$ is surjective there are elements in $M \circ N$ which are represented as $[m,n]$ and $[m,n']$ for some $m \in M$. Since $l_{M \circ N}([m,n]) = l_{M}(m)$ the elements $[m,n]$ and $[m,n']$ are in the same left fiber of $M \circ N$. By assumption the right $G$ action on $M \circ N$ is transitive, so there is a $g \in G$ such that $[m,n] \cdot g = [m,n']$. This implies that there is an $h \in H$ such that $(m, n \cdot g) = (m \cdot h^{-1}, h \cdot n')$. But since the right action of $H$ on $M$ is free, $m = m \cdot h^{-1}$ implies $h = 1_{r_{M}(m)} = 1_{l_{N}(n')}$ and, hence, $h \cdot n' = n'$. Thus $n \cdot g = n'$. We conclude that $G$ acts transitively on the left fibers of $N$. By analogous arguments we can show that $N$ and $M$ are left and right transitive.
\end{proof}

\begin{Corollary}
\label{th:isoprincipal2}
Let $M$ be a smooth bibundle. If $N$ is a weak inverse of $M$ then $N \cong M^{\op}$ 
\end{Corollary}

\section{Stacky Lie groups}
\label{sec:Groups}

We want to gain insight into the constructive aspects of differentiable stacks by studying their presentations by Lie groups. In this spirit we could define a ``stacky Lie group'' to be the presentation of a differentiable group stack \cite{Breen:Bitorseurs}. But this definition isn't of much value if we want to study such presentations in their own right. What we need is an axiomatic definition of a stacky Lie group in terms of a Lie groupoid and group structure morphisms (i.e. bibundles and isomorphisms of bibundles) subject to certain relations.

A differentiable group stack is a group object internal to the category of differentiable stacks. As we have seen, the category of differentiable stacks extends to a (strict) 2-category which is equivalent to the weak 2-category $\stMfld$. By this equivalence a group stack is mapped to a Lie groupoid which is equipped with a group structure in some sense. More precisely, we should get a group object internal to the weak 2-category $\stMfld$, a notion which is yet to be defined.

There is a considerable amount of literature on various notions of weak and higher groups. However, to our best knowledge the only concept of group objects \emph{in} weak higher categories has been given implicitly in Lurie's generalized definition of groupoids of $\infty$-categories \cite{Lurie:Topos}. His definition is too general for our purpose, however, because it works best in categories with finite limits. In the category of manifolds it does not reproduce the usual notion of Lie groupoid. Our definition is a modification of his definition. For the special case of 1-group objects in a 2-category, which is relevant here, we also spell out the simplicial construction and give a definition in terms of structure morphisms and coherence relations.

In any case, the concept of groups in higher categories is somewhat involved and quite an interesting subject in its own right. In this paper we can only give a concise exposition of the general construction. A full account of the general theory will be given elsewhere \cite{Blohmann:Higher}.

\subsection{Simplicial objects in higher categories}

Recall, that a simplicial set is a functor $X: \sDelta^\op \rightarrow \Set$ where $\sDelta$ denotes the category of finite ordered sets and weakly monotonous maps. The morphisms between two simplicial sets $X$ and $Y$ are the natural transformations from $X$ to $Y$. The objects in $\sDelta$ are denoted by $[n]= \{0,1,2,\ldots,n\}$. The image of $[n]$ under $X$ is $X_n := X([n])$. The simplicial set represented by $[n]$ is denoted by $\Delta^n := \Hom_{\sDelta}(\Empty, [n])$. Given a simplicial set $X$ the Yoneda lemma identifies $\Hom_{\sSet}(\Delta^n,X) \cong X_n$.

To every small category $\catC$ we can associate a simplicial set $\Nerve(\catC)$, called the nerve of $\catC$, given by mapping $[0]$ to the set of objects of $\catC$ and $[n]$, $n > 0$ to the set of sequences of $n$ composable morphisms $C_0 \rightarrow C_1 \rightarrow \ldots \rightarrow C_n$ between objects in $\catC$. In order to decide whether a simplicial set is the nerve of some category we can use the concept of Kan conditions.

The representable simplicial set $\Delta^n$ has simplicial subsets $\Lambda^n_i$, $0 \leq i \leq n$ which are generated by removing the interior and the $i$-th face. $\Lambda^n_i$ is called the $i$-th horn. $\Lambda^n_0$ and $\Lambda^n_n$ are called outer horns, all other horns are called inner. There is an obvious inclusion of simplicial sets $\Lambda^n_i \hookrightarrow \Delta^n$. We say that a simplicial set $X$ satisfies the weak Kan condition $\Kan(n,i)$ if the induced map $\Hom_{\sSet}(\Lambda^n_i, X) \rightarrow \Hom_{\sSet}(\Delta^n, X) \cong X_n$ is surjective. We say that $X$ satisfies the strict Kan condition $\Kan!(n,i)$ if the map is a bijection.

Kan conditions can be used to describe simplicial sets which are the nerves of categories or groupoids, i.e. classifying spaces. In fact, a simplicial set is the nerve of a category iff it satisfies the strict inner Kan condition $\Kan!(k,i)$, $0 < i < k$, for all $k \geq 2$. This generalizes naturally as follows: A Joyal $n$-category (also called quasi-category) is a simplicial set $\catC$ which satisfies the weak inner Kan condition $\Kan(k,i)$, $0 < i < k$, for $k \leq n$ and the strict inner Kan conditions $\Kan!(k,i)$, $0 < i < k$ for $k \geq n$. In this sense ordinary categories can be identified with Joyal 1-categories. A Joyal 2-category is the the same as a bicategory (in the sense of B\'enabou) in which all 2-morphisms are isomorphisms. 
In the following ``$n$-category'' will always mean ``Joyal $n$-category''.

Kan conditions can also be used to describe set-theoretic groupoids as simplicial sets which satisfy the strict Kan condition $\Kan!(k,i)$, $0 \leq i \leq k$, for all $k \geq 2$. The outer Kan conditions guarantee the existence of inverses. This generalizes to $n$-groupoids which are defined to be simplicial sets which satisfy the weak Kan condition $\Kan(k,i)$, $0 \leq i \leq k$ for $2 \leq k \leq n$ and the strict Kan conditions $\Kan!(k,i)$, $0 \leq i \leq k$ for $k > n$. A set theoretic $n$-group, finally, is an $n$-groupoid $G$ with a single $0$-simplex, $G_0 = *$.

In order to internalize this simplicial description of groups to higher categories we have two answer two questions: What is a simplicial object internal to an $n$-category and what are the Kan conditions for such a simplicial object? To answer the first question, we observe that we can identify a functor of categories with a morphism of simplicial sets of the nerves of the categories. In fact, for any 1-category $\calC$ there is a natural bijection
$\Hom_{\Cat}(\sDelta^\op, \catC) \cong \Hom_{\sSet}(\Nerve(\sDelta^\op), \Nerve(\catC))$. This leads to the following definition:

\begin{Definition}
  A \emph{simplicial object in an $n$-category} $\catC$ is a morphism of simplicial sets $\Nerve(\sDelta^\op) \rightarrow \catC$.
\end{Definition}

This is also the definition used by Lurie (definition~6.1.2.2 in \cite{Lurie:Topos}). Explicitly such a higher simplicial object $G$ maps a sequence of $k$ composable morphisms in $\Delta^\op$ to a $k$-simplex in $\catC$,
\begin{equation*}
  \Nerve(\sDelta)_k \ni
  \bigl( 
  [n_0] \stackrel{f_1}{\leftarrow}
  [n_1] \stackrel{f_2}{\leftarrow}
  [n_2] \stackrel{f_2}{\leftarrow} \ldots
  \stackrel{f_k}{\leftarrow} [n_k] \bigr)
  \longmapsto 
  G(f_1,f_2,\ldots,f_k) \in \catC_k \,.
\end{equation*}
The images of the identity morphisms will be denoted by $G(\id_{[n]}) \equiv G_n \in \catC_0$. If $\catC = \Nerve(\catD)$ a 1-category we retrieve the usual notion of simplicial objects in $\catD$.

\subsection{Groups in higher categories}
\label{sec:GroupsHigher}

The question of how to describe the Kan conditions for a simplicial object in a higher category is more intricate. Already for a simplicial object $X: \sDelta^\op \rightarrow \catC$ in a category $\catC$ it no longer makes sense to define the horns as hom-space ``$\Hom(\Lambda^n_i,X)$'' because $\Lambda^n_i$ and $X$ are functors to different categories. Even if there is a forgetful functor from $\catC$ to $\Set$ we have to enrich the set ``$\Hom_{\sSet}(\Lambda^n_i,X)$'' over $\catC$. 

A good way to define the horn spaces is as limits in $\catC$. Toward this end we recall the definition of the simplex category $(\sDelta \Comma S)$ of a given simplicial set $S$: The objects are morphisms of simplicial sets $\Delta^k \rightarrow S$, the morphisms are commutative triangles
\begin{equation*}
\xymatrix@-1ex@C-2ex{
  \Delta^k \ar[rr]^{\phi} \ar[dr] && \Delta^{k'} \ar[dl] \\
 & S} 
\end{equation*}
The Yoneda lemma identifies the morphism $\phi$ with a morphism $\phi:[k] \rightarrow [k']$ in $\sDelta$, so we have a natural projection $(\sDelta \Comma S) \rightarrow \sDelta$. If $G$ is another simplicial set this induces a functor
\begin{equation*}
  G[S] : (\sDelta \Comma S)^\op
 \longrightarrow \sDelta^\op
 \stackrel{G}{\longrightarrow} \Set \,,
\end{equation*}
using the notation of Lurie \cite{Lurie:Topos}. The set of morphisms from $S$ to $G$ is the limit of this functor,
\begin{equation*}
  \Hom_{\sSet}(S,G) \cong \varprojlim G[S] \,.
\end{equation*}
This definition generalizes to simplicial objects in other categories than $\Set$. For example for the space of inner 2-horns the limit, if it exists, is given by
\begin{equation*}
  \varprojlim G[\Lambda^2_1] = G_1 \times_{G_0}^{d_0,d_1} G_1 \,.
\end{equation*}
The definition of horn-spaces as limits generalizes to simplicial objects in higher categories. First, we need to recall the notion of limit in a higher category as terminal cone, which is due to Joyal \cite{Joyal:Quasi} (cf. also Lurie \cite{Lurie:Topos}). 

On $\sDelta$ we have the ordinal sum $[m] + [n] = [m+n+1]$ which extends naturally by the Yoneda embedding $\Delta^m \star \Delta^n = \Delta^{m+n+1}$ to an operation $\star: \sSet \times \sSet \rightarrow \sSet$ given by 
\begin{equation*}
  (S \star T)_n = S_n \sqcup T_n \sqcup \bigsqcup_{i+j = n -1} S_i \times T_j \,,
\end{equation*}
for $S,T \in \sSet$. This operation is called the \emph{join} because it is the simplicial analogue of the topological join. Explicitly, the join of $S$ and $T$ consist of all simplices of $S$, all simplices of $T$, all edges joining vertices of $S$ to vertices of $T$, all 2-simplices with edges joining $S$ and $T$ etc. The join $\Delta^0 \star S$ of a point with a simplicial set can be pictured as cone with tip $\Delta^0 = *$, basis $S$, and edges going from the tip to all vertices of $S$. (Lurie calls this this the left cone, Joyal the projective cone.) 

Given a morphism $D:S \rightarrow T$ of simplicial sets, viewed as $S$-parametrized diagram in $T$, we can define a simplicial set $T_{/D}$ by
\begin{equation*}
  (T_{/D})_n := \Hom^D_{\sSet}(\Delta^n \star S, T) := 
  \{F: \Delta^n \star S \rightarrow T \,:\, F|_{S} = D \} \,.
\end{equation*}
It can be shown (\cite{Lurie:Topos}, prop.~1.2.9.3) that if $T$ is an $\infty$-category then $S_{/D}$ is also an $\infty$-category, called the \emph{overcategory} of $D$. If $S$ and $T$ are the nerves of ordinary categories, then $T_{/D}$ is the nerve of the overcategory in the usual sense.

In an ordinary category a limit of $D$ is a terminal object in the overcategory. An object $c \in \catC_0$ of an $\infty$-category $\catC$ is called terminal if every simplicial map $F: \partial \Delta^n \rightarrow \catC$ with $F(0) = c$ has an extension to a simplex $\tilde{F}: \Delta^n \rightarrow \catC$. Now we can finally give the definition we need:

\begin{Definition}
  Let $S$ be a simplicial set, $\catC$ an $\infty$-category, and $D: S \rightarrow \catC$ a map of simplicial sets. A \emph{limit} of $D$ is a terminal object in $\catC_{/D}$.
\end{Definition}

We denote the limit by $\varprojlim D$, which is a customary abuse of language since the limit is not unique and involves a choice. But it can be shown that in $\infty$-categories the limits of a given diagram generate a Kan complex which is contractible, i.e., homotopic to a point. The notion of $\infty$-limit specializes in an obvious way to $n$-categories. 

For a simplicial object $G$ in an $\infty$-category $\catC$ and a simplicial set $S$ there is a map of simplicial sets $G[S] : \Nerve(\sDelta \Comma S)^\op \rightarrow \catC$ induced by the projection of the simplex category of $S$. If the limit of this map exists it will be denoted by 
\begin{equation*}
  G^S := \varprojlim G[S] \,.
\end{equation*}
If $S$ is finite and $\catC$ has finite coproducts this limit can be understood as right exponential object which justifies the notation. Now we can define Kan conditions in higher categories.

\begin{Definition}
  Let $G: \Nerve(\sDelta)^\op \rightarrow \catC$ be a simplicial object in an $\infty$-category $\catC$. If the limit of $G[\Lambda^n_i]$ exists and if the morphism 
\begin{equation*}
  G_n \longrightarrow G^{\Lambda^n_i} \,,
\end{equation*}
induced by the inclusion $\Lambda^n_i \hookrightarrow \Delta^n$, has a section then we say that $G$ satisfies the weak Kan condition $\Kan(n,i)$. If it is an isomorphism we say that $G$ satisfies the strict Kan condition $\Kan!(n,i)$
\end{Definition}

There are several ways to further relax these notions of Kan conditions. For the weak Kan condition we could alternatively require the morphism to be epi. Our definition has the advantage that it is preserved by any functor of the target category $\catC$ which preserves finite limits. Lurie does not require the existence of the horn-space limits and only requires the map $G[\Delta^n] \rightarrow G[\Lambda^n_i]$ to be epi. For our purposes the definition given here is more appropriate because it reproduces the usual definition of Lie groupoid in $\Mfld$. In any case the definition of groupoids in terms of Kan conditions takes on the  usual form:

\begin{Definition}
  Let $G: \Nerve(\sDelta^\op) \rightarrow \catC$ be a simplicial object internal to the $\infty$-category $\catC$. If $G$ satisfies the weak Kan condition $\Kan(k,i)$, $0 \leq i \leq k$ for all $k \leq n$ and the strict Kan condition $\Kan!(k,i)$, $0 \leq i \leq k$ for all  $k > n$ then we call $G$ an $n$-groupoid in $\catC$. If furthermore $G_0$ is terminal in $\catC$, then we call $G$ an $n$-group in $\catC$.
\end{Definition}

If $\catC$ is a $1$-category this definition reproduces the usual definition of $n$-groups and $n$-groupoids. If we require only the inner Kan conditions to hold we obtain the definition of an $n$-category and an $n$-monoid in $\catC$ in complete analogy to the usual case. The following proposition is a simple observation.

\begin{Proposition}
\label{th:FunctGroup}
  Let $G$ be an $n$-group in the $\infty$-category $\catC$ and $f: \catC \rightarrow \catD$ a functor of $\infty$-categories which preserves finite limits. Then $f G$ is an $n$-group in $\catD$.
\end{Proposition}

This statement holds in particular if $f$ is an equivalence of higher categories.

\subsection{Stacky Lie groups}

Note that all the statements of the preceding section also make sense if $\catC$ is an $n$-category. Now we can finally give the sought-after definition.

\begin{Definition}
  A \emph{stacky Lie group} is a 1-group object in the 2-category $\stMfld$.
\end{Definition}

By the 2-categorical equivalence of $\stMfld$ and differentiable stacks we obtain the following corollary to proposition~\ref{th:FunctGroup}:

\begin{Corollary}
\label{th:Coro}
  Every presentation of a differentiable group stack is a stacky Lie group.
\end{Corollary}

As promised, we will now spell out this definition in terms of structure morphisms and coherence relations. Let $G$ be a stacky Lie group. Since $\stMfld$ is a 2-category it is 2-coskeletal. This means that the morphism $G$ is determined by its 2-truncation, i.e., by its values on the 0, 1, and 2-simplices of $\Nerve(\sDelta)^\op$. The image of an object is a Lie groupoid $G_n = G([n])$, the image of a morphism $\phi: [n] \rightarrow [n']$ a right principal bibundle. The image of a pair of composable arrows $\phi$ and $\psi$ can be identified with an isomorphism of bibundles $G(\phi,\psi): G(\phi) \circ G(\psi) \rightarrow G(\phi\psi)$. (Recall, that we compose bibundles from left to right.) The set of morphisms of $\sDelta$ is generated by the coface maps $d^i$ and the codegeneracy maps $s^j$. As usual, we denote their images by $d_i := G(d^i)$ and $s_j := G(s^j)$. Every bibundle in the image of $G$ is isomorphic to a product of these morphisms.

The space of inner 2-horns can be identified with $G^{\Lambda^2_1} = G_1 \times_{G_0} G_1 = G_1 \times G_1$. More generally, it can be shown by induction that we can identify $\Lambda^k_i(G) = (G_1)^{k}$. The Kan condition then tells us that $G_k \cong (G_1)^{k}$. The face morphism
\begin{equation*}
  \mu: G_1 \times G_1 = G^{\Lambda^2_1} 
  \stackrel{\cong}{\longrightarrow} G_2 
  \stackrel{d_1}{\longrightarrow} G_1 
\end{equation*}
is the structure morphism of multiplication. The group identity is given by the degeneracy morphism
\begin{equation*}
  e: 1 \cong G_0 \stackrel{s_0}{\longrightarrow} G_1 \,.
\end{equation*}
So far everything looks like for a group object in an ordinary category. However, even though we have a 1-group object we do not have strict associativity for the group multiplication. In fact, from the identity $d^2 d^1 = d^1 d^1$ of coface maps it merely follows that there is a 2-isomorphism
\begin{equation}
\label{eq:preass}
  d_2 \circ d_1 
  = G(d^2) \circ G(d^1) 
  \stackrel{\cong}{\longrightarrow}
  G(d^2 d^1) \stackrel{\cong}{\longrightarrow}
  G(d^1 d^1) \stackrel{\cong}{\longrightarrow}
  G(d^1) \circ G(d^1) = 
  d_1 \circ d_1 \,.
\end{equation}
To see how we get from this relation the group associator we identify as before $G_k \cong (G_1)^{k}$. The inner face maps $d^{[k]}_i: (G_1)^{k} \rightarrow (G_1)^{k-1}$ can then be identified with
\begin{equation}
\label{eq:faceabbr}
\begin{gathered}
  d^{[2]}_1 \cong \mu \,, \\
  d^{[3]}_1 \cong \mu \times \Id \,,\quad 
  d^{[3]}_2 \cong \Id \times \mu \,, \\
  d^{[4]}_1 \cong \mu \times \Id \times \Id \,,\quad 
  d^{[4]}_2 \cong \Id \times \mu \times \Id \,,\quad 
  d^{[4]}_3 \cong \Id \times \Id \times \mu \,,
\end{gathered}
\end{equation}
and so forth. The outer face maps $d^{[k]}_0$ and $d^{[k]}_k$ are the projections on the first $k-1$ and the last $k-1$ factors of $(G_1)^k$, respectively. Similarly, the degeneracy maps can be identified with
\begin{equation}
\label{eq:degabbr}
\begin{gathered}
  s^{[0]}_0 \cong e \,, \\
  s^{[1]}_0 \cong e \times \Id \,,\quad 
  s^{[1]}_1 \cong \Id \times e \,, \\
  s^{[2]}_0 \cong e \times \Id \times \Id \,,\quad 
  s^{[2]}_1 \cong \Id \times e \times \Id \,,\quad 
  s^{[2]}_2 \cong \Id \times \Id \times e \,,
\end{gathered}
\end{equation}
and so forth.

As in \eqref{eq:preass} we obtain as images of various simplicial relations the 2-isomorphisms of the associator, the left and the right unit constraints:
\begin{equation}
\label{eq:stackyass}
\begin{aligned}
  \ass : 
  (\mu \times \Id) \circ \mu \cong 
  d^{[3]}_1 \circ d^{[2]}_1 
  &\stackrel{\cong}{\longrightarrow}
  d^{[3]}_2 \circ d^{[2]}_1 \cong
  (\Id \times \mu) \circ \mu \,,\\
  \luni : 
  (e \times \Id) \circ \mu \cong 
  s^{[1]}_0 \circ d^{[2]}_1 
  &\stackrel{\cong}{\longrightarrow}
  \Id \,,\\
  \runi : 
  (\Id \times e) \circ \mu \cong 
  s^{[1]}_1 \circ d^{[2]}_1 
  &\stackrel{\cong}{\longrightarrow}
  \Id \,.\\
\end{aligned}
\end{equation}
The other simplicial relations yield 2-isomorphisms such as 
\begin{equation}
\label{eq:MissingIso}
  (\mu \times \Id \times \Id) \circ (\Id \times \mu)
  \cong d_1^{[4]} \circ d_2^{[3]} 
  \stackrel{\cong}{\longrightarrow}
  d_3^{[4]} \circ d_1^{[3]} \cong
  (\Id \times \Id \times \mu) \circ (\mu \times \Id) \,,
\end{equation}
which are implicit in the identifications \eqref{eq:faceabbr} and \eqref{eq:degabbr} and are, therefore, no further data of the stacky Lie group structure. Just as for monoids in an ordinary category, the fundamental coherence relation comes from the Kan condition for filling a 4-horn: 
\begin{equation}
\label{eq:weakpentagon}
\begin{gathered}
\xymatrix@C+4ex{
 (d_1^{[4]} \circ d_1^{[3]}) \circ d_1^{[2]} \ar[r]^{\Assc} \ar[d]_{(\ass \times \id) \circ \id~} & d_1^{[4]} \circ (d_1^{[3]} \circ d_1^{[2]}) \ar[d]^{\id \circ \ass}  
\\
(d_2^{[4]} \circ d_1^{[3]}) \circ d_1^{[2]} \ar[d]_{\Assc} & d_1^{[4]} \circ (d_2^{[3]} \circ d_1^{[2]}) \ar[d]^{\Assc^{-1}}
\\
d_2^{[4]} \circ (d_1^{[3]} \circ d_1^{[2]}) \ar[d]_{\id \circ \ass} & (d_1^{[4]} \circ d_2^{[3]}) \circ d_1^{[2]} \ar[d]^{\cong} 
\\
d_2^{[4]} \circ (d_2^{[3]} \circ d_1^{[2]}) \ar[d]_{\Assc^{-1}} & (d_3^{[4]} \circ d_1^{[3]}) \circ d_1^{[2]} \ar[d]^{\Assc} 
\\
(d_2^{[4]} \circ d_2^{[3]}) \circ d_1^{[2]} \ar[d]_{~~(\id \times \ass) \circ \id} & d_3^{[4]} \circ (d_1^{[3]} \circ d_1^{[2]}) \ar[d]^{\id \circ \ass} 
\\
(d_3^{[4]} \circ d_2^{[3]}) \circ d_1^{[2]} \ar[r]^{\Assc} & d_3^{[4]} \circ (d_2^{[3]} \circ d_1^{[2]})
}
\end{gathered}
\end{equation}
Here, the associator $\ass$ of the group is distinguished by a subscript from the associator $\Assc$ of the 2-category $\stMfld$. The unlabeled isomorphism on the right hand side is \eqref{eq:MissingIso}. Note, that if the composition of bibundles were strictly associative, this dodecagon would reduce to the usual pentagon. The relation for the unit constraints is:
\begin{equation}
\label{eq:weakunit}
\begin{gathered}
\xymatrix@C-0ex{
  s_1^{[2]} \circ (d^{[3]}_1 \circ d^{[2]}_1)
  \ar[rr]^{\id \circ \ass}
  \ar[d]_{\Assc}
& &
  s_1^{[2]} \circ (d^{[3]}_2 \circ d^{[2]}_1)
  \ar[d]^{\Assc}
\\
  (s_1^{[2]} \circ d^{[3]}_1) \circ d^{[2]}_1
  \ar[dr]_{(\runi \times \id) \circ \id~~}
& &
  (s_1^{[2]} \circ d^{[3]}_2) \circ d^{[2]}_1
  \ar[dl]^{~~(\id \times \luni) \circ \id}
\\
& d^{[2]}_1 &
}
\end{gathered}
\end{equation}
This pentagon reduces to the usual triangle relation if the composition of bibundles is associative on the nose. There are more such coherence diagrams which can be deduced from \eqref{eq:weakpentagon} and \eqref{eq:weakunit} together with the relations implicit in the special form \eqref{eq:faceabbr} and \eqref{eq:degabbr} of the structure morphisms. These notably include the diagrams corresponding to the other two triangle diagrams and a diagram corresponding to the equality of the left and right unit constraints on the square of the group identity.

Until now, we have used the Kan condition for inner horns in order to identify $G_2$ with cartesian powers of $G_1$. This corresponds to the bar resolution of a semi-group. It remains to explain how the structure morphism of the inverse arises from the outer horn-filling conditions. By definition, we have isomorphisms $\sigma_i : G^{\Lambda^2_i} \rightarrow G_2$ for $i=0,1,2$ which fit into the following diagram in $\stMfld$:
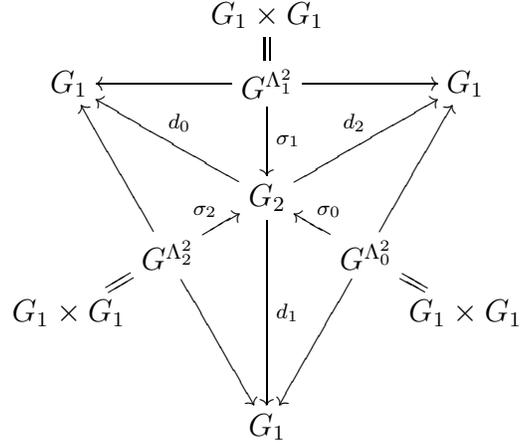
\begin{equation*}
\xy 0;<15pt,0pt>:;
  (-5,0)*+{G_1}="D0";
  (+5,0)*+{G_1}="D2";
  (0,-8.66)*+{G_1}="D1";
  "D1";"D2"**{}; ?(0.5)*+{G^{\Lambda^2_0}}="L0";
  "D0";"D2"**{}; ?(0.5)*+{G^{\Lambda^2_1}}="L1";
  "D0";"D1"**{}; ?(0.5)*+{G^{\Lambda^2_2}}="L2";
  "D0";"L0"**\dir{} ?! {"D1";"L1"} *+{G_2}="C";
  "C";"L0"**{}; ?(2)*+{G_1\times G_1}="G0";
  "C";"L1"**{}; ?(1.6)*+{G_1\times G_1}="G1";
  "C";"L2"**{}; ?(2)*+{G_1\times G_1}="G2";
  \ar @{->} "L1"; "D0"
  \ar @{->} "L1"; "D2"
  \ar @{->} "L2";"D0"
  \ar @{->} "L2";"D1"
  \ar @{->} "L0";"D1"
  \ar @{->} "L0";"D2"
  \ar @{->}_{d_0} "C";"D0"
  \ar @{->}^{d_1} "C";"D1"
  \ar @{->}^{d_2} "C";"D2"
  \ar @{<-}^{\sigma_0} "C";"L0"
  \ar @{<-}_{\sigma_1} "C";"L1"
  \ar @{<-}_{\sigma_2} "C";"L2"
  \ar @{=} "G0";"L0"
  \ar @{=} "G1";"L1"
  \ar @{=} "G2";"L2"
\endxy
\end{equation*}
Every edge is a bibundle and every triangle represents a isomorphism of bibundles. We have omitted the 2-isomorphisms in order to improve legibility. From this 2-commutative diagram we get to the usual description of a group by the structure morphisms of the left and right inverse after a diagram chase:
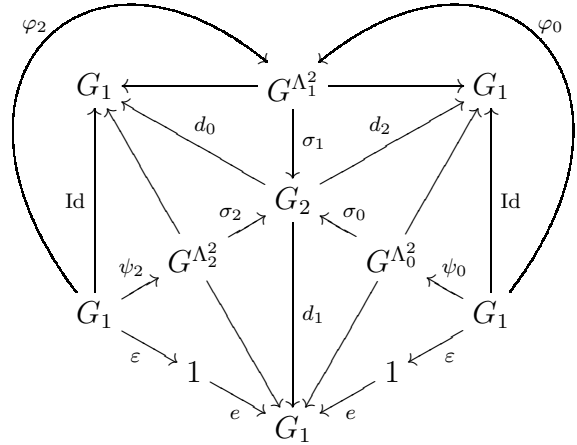
\begin{equation*}
\xy 0;<15pt,0pt>:;
  (-5,0)*+{G_1}="D0";
  (+5,0)*+{G_1}="D2";
  (0,-8.66)*+{G_1}="D1";
  "D1";"D2"**{}; ?(0.5)*+{G^{\Lambda^2_0}}="L0";
  "D0";"D2"**{}; ?(0.5)*+{G^{\Lambda^2_1}}="L1";
  "D0";"D1"**{}; ?(0.5)*+{G^{\Lambda^2_2}}="L2";
  "D0";"L0"**\dir{} ?! {"D1";"L1"} *+{G_2}="C";
  "C";"L0"**{}; ?(2)*+{G_1}="G0";
  "G0";"D1"**{}; ?<>(0.5)*+{1}="1";
  "G0";"L1" **\crv{(10,0)&(5,5)} ?>*\dir{>};
  (6.5,1.5)*+{\scriptstyle\phi_0};
  "C";"L2"**{}; ?(2)*+{G_1}="G2";
  "G2";"D1"**{}; ?<>(0.5)*+{1}="2";
  "G2";"L1" **\crv{(-10,0)&(-5,5)} ?>*\dir{>};
  (-6.5,1.5)*+{\scriptstyle\phi_2};
  \ar @{->} "L1"; "D0"
  \ar @{->} "L1"; "D2"
  \ar @{->} "L2";"D0"
  \ar @{->} "L2";"D1"
  \ar @{->} "L0";"D1"
  \ar @{->} "L0";"D2"
  \ar @{->}_{d_0} "C";"D0"
  \ar @{->}^{d_1} "C";"D1"
  \ar @{->}^{d_2} "C";"D2"
  \ar @{<-}^{\sigma_0} "C";"L0"
  \ar @{<-}_{\sigma_1} "C";"L1"
  \ar @{<-}_{\sigma_2} "C";"L2"
  \ar @{->}^{\varepsilon} "G0";"1"
  \ar @{->}^{e} "1";"D1"
  \ar @{->}_{\Id} "G0";"D2"
  \ar @{->}_{\psi_0} "G0";"L0"
  \ar @{->}_{\varepsilon} "G2";"2"
  \ar @{->}_{e} "2";"D1"
  \ar @{->}^{\Id} "G2";"D0"
  \ar @{->}^{\psi_2} "G2";"L2"
\endxy
\end{equation*}
Here, $\varepsilon$ is the terminal morphism. The morphisms $\psi_0$ and $\psi_2$ are given by the universal property of the limits $G^{\Lambda^2_0}$ and $G^{\Lambda^2_2}$. The structure bibundles of the left and right inverse are given by the composition
\begin{equation*}
\begin{aligned}
  i_l := (\psi_0 \circ \sigma_0) \circ d_0 \,,\qquad
  i_r := (\psi_2 \circ \sigma_2) \circ d_2 \,.
\end{aligned}
\end{equation*}
The morphisms $i_l$ and $\Id$ induce by the universal property of $G^{\Lambda^2_1}$ the morphism $\phi_0 = \Delta \circ (i_l \times \Id)$. Analogously, we obtain $\phi_2 = \Delta \circ (\Id \times i_r)$. From the diagram we read off the 2-commutative diagram
\begin{equation}
\label{eq:invcoherence}
\begin{gathered}
  \xymatrix@dr@+4ex{
  G_1 
  \ar[dr]_(0.65){\varepsilon \circ e} 
  \ar[d]_{\Delta \circ (i_l \times \Id)} 
  \ar[r]^{\Delta \circ (\Id \times i_r)} 
  \drtwocell<\omit>{^<3>{\linv}}
  \drtwocell<\omit>{_<-3>{\rinv}}
  & G_1 \times G_1 \ar[d]^{\mu}\\
  G_1 \times G_1 \ar[r]_{\mu} & G_1
}
\end{gathered}
\end{equation}
The 2-isomorphisms of this diagram are called the left and right inverse constraints. It is shown below that the existence of a left inverse implies the existence of a right inverse to which it is isomorphic. We wrap up the description of a stacky Lie group by structure morphisms:

\begin{Proposition}
  The following data determine the structure of a stacky Lie group: A Lie groupoid $G$;  a $(G\times G)$-$G$ bibundle $\mu$, a $1$-$G$ bibundle $e$, and two $G$-$G$ bibundles $i_l$ and $i_r$; isomorphisms of bibundles
\begin{equation*}
\begin{gathered}
  \ass : (\mu \times \Id) \circ \mu 
  \stackrel{\cong}{\longrightarrow} (\Id \times \mu) \circ \mu \,,\\
  \luni : (e \times \Id) \circ \mu 
  \stackrel{\cong}{\longrightarrow} \Id \,,\qquad
  \runi : (\Id \times e) \circ \mu
  \stackrel{\cong}{\longrightarrow} \Id \,,\\
  \linv: (\Delta \circ (i_l \times \id)) \circ \mu 
  \stackrel{\cong}{\longrightarrow} \varepsilon \circ e \,,\qquad
  \rinv: (\Delta \circ (\id \times i_r)) \circ \mu
  \stackrel{\cong}{\longrightarrow} \varepsilon \circ e \,,
\end{gathered}
\end{equation*}
which satisfy the coherence relations \eqref{eq:weakpentagon}
 and \eqref{eq:weakunit}, using the abbreviations \eqref{eq:faceabbr} and \eqref{eq:degabbr}.
\end{Proposition}

\subsection{The graphic notation}

For many interesting questions the 2-morphism do not matter and we can work in the 1-category of Lie groupoids and isomorphism classes of principal bibundles. In the graphic notation of the previous section we represent the isomorphism classes of the structure bibundles of a stacky group by:
\begin{alignat*}{2}
&\text{multiplication:}& [\mu] &=
\begin{xyc}
  \Copr(0,0)
\end{xyc}\\[1ex]
&\text{neutral element:}\quad& [e] &=
\begin{xyc}
  \Coun(0,0)
\end{xyc}\\[1ex]
&\text{inverse:}& [i] &=
\begin{xyc}
  \Anti(0,0)
\end{xyc}
\end{alignat*}
In this notation the commutative diagrams of associativity and unitality read
\begin{equation}
\label{eq:monoidal}
\begin{xyc}
  \Copr(0,0) \Copr(-1,2) \Iden(1,1)\Iden(1,2)
\end{xyc}
\quad = \quad
\begin{xyc}
  \Copr(0,0) \Copr(1,2) \Iden(-1,1)\Iden(-1,2)
\end{xyc}
\qquad\text{and}\qquad
\begin{xyc}
  \Copr(0,0) \Coun(-1,1) \Iden(1,1)
\end{xyc}
\quad = \quad
\begin{xyc}
  \Iden(0,-1)\Iden(0,0)\Iden(0,1)
\end{xyc}
\quad = \quad
\begin{xyc}
  \Copr(0,0) \Coun(1,1) \Iden(-1,1)
\end{xyc} \quad.
\end{equation}
The diagram of the inverse becomes
\begin{equation}
\label{eq:antipode}
\begin{xyc}
  \Copr(0,0) \Anti(-1,1) \Iden(1,1) \Prod(0,3)
\end{xyc}
\quad = \quad
\begin{xyc}
  \Coun(0,0) \Unit(0,3)
\end{xyc}
\quad = \quad
\begin{xyc}
  \Copr(0,0) \Iden(-1,1) \Anti(1,1) \Prod(0,3)
\end{xyc} \quad.
\end{equation}
Moreover, we can spell out that the structure bibundles are right principal using Prop.~\ref{th:rightprincipal}. In graphic notation this reads for the multiplication:
\begin{equation}
\label{eq:bimonoid1}
\begin{xyc}
 \Prod(0,1)\Coun(0,1)
\end{xyc}
\quad = \quad
\begin{xyc}
 \Iden(-1,0)\Coun(-1,1) \Iden(1,0)\Coun(1,1)
\end{xyc}
\qquad\text{and}\qquad
\begin{xyc}
 \Copr(-2,0)\Copr(2,0)
 \Iden(-3,1)\Iden(-3,2) \Flip(0,2) \Iden(3,1)\Iden(3,2)
 \Prod(-2,4)\Prod(2,4)
\end{xyc}
\quad = \quad
\begin{xyc}
 \Prod(0,0)\Copr(0,2) 
\end{xyc} \quad,
\end{equation}
the neutral element:
\begin{equation}
\label{eq:bimonoid2}
\begin{xyc}
 \Term(0,1)\Copr(0,2)
\end{xyc}
\quad = \quad
\begin{xyc}
 \Unit(-1,0)\Iden(-1,1) \Unit(1,0)\Iden(1,1)
\end{xyc}
\qquad\text{and}\qquad
\begin{xyc}
 \Unit(0,0)\Coun(0,1)
\end{xyc}
\quad = \quad
\begin{xyc}
  \Term(0,0)
\end{xyc}
\quad, 
\end{equation}
and the inverse:
\begin{equation}
\label{eq:antipodecounital}
\begin{xyc}
  \Unit(0,0)\Anti(0,1)
\end{xyc}
\quad = \quad
\begin{xyc}
  \Unit(0,0)\Iden(0,1)
\end{xyc}
\qquad\text{and}\quad
\begin{xyc}
  \Prod(0,1)\Anti(0,2)
\end{xyc}
\quad = \quad
\begin{xyc}
  \Anti(-1,0)\Anti(1,0)\Prod(0,2)
\end{xyc}
\quad.
\end{equation}

These diagrams can also be interpreted in the following way: Diagrams~\eqref{eq:monoidal} mean that $\mu$ and $e$ equip $G$ with a weak monoid structure. Diagrams~\eqref{eq:comonoidal} state that $\Delta$ and $\varepsilon$ equip $G$ with a weak comonoid structure. Diagrams~\eqref{eq:bimonoid1} and \eqref{eq:bimonoid2} tell us that the monoid structure and the comonoid structure are compatible. More precisely, $\Delta$ and $\varepsilon$ are monoidal morphisms whereas $\mu$ and $e$ are comonoidal. This can be expressed by saying that $G$ is a bimonoid internal to the category of Lie groupoids and isomorphism classes of principal bibundles.

\subsection{The preinverse}

A bialgebra is a bimonoid in the category of vector spaces. In analogy to the antipode of a Hopf algebra we then expect the inverse of the stacky group to be a unital and counital but antimonoidal and anticomonoidal morphism. Indeed, counitality is expressed by the first of diagrams~\eqref{eq:antipodecounital}. Unitality is shown by the following graphic calculation:
\begin{equation*}
\begin{xyc}
  \Anti(0,0) \Coun(0,1)
\end{xyc}
\quad = \quad
\begin{xyc}
  \Coun(0,1) \Coun(2,1)
  \Anti(0,0) \Iden(2,0)
  \Copr(1,-1)
\end{xyc}
\quad = \quad
\begin{xyc}
  \Coun(0,3)
  \Copr(0,0) \Anti(-1,1) \Iden(1,1) \Prod(0,3)
\end{xyc}
\quad = \quad
\begin{xyc}
  \Coun(0,0) \Unit(0,3) \Coun(0,4)
\end{xyc}
\quad = \quad
\begin{xyc}
  \Iden(0,0) \Coun(0,1)
\end{xyc}
\quad,
\end{equation*}
where we have used \eqref{eq:monoidal}, \eqref{eq:bimonoid1}, \eqref{eq:antipode}, and \eqref{eq:bimonoid2}. The second diagram of \eqref{eq:antipodecounital} states that the inverse is a weakly comonoidal morphism. But since the comonoid product $\Delta$ is cocommutative a weakly comonoidal morphism is the same thing as a weakly anticomonoidal morphism. The stacky inverse is also a weakly antimonoidal morphism as the following lemma shows:

\begin{Lemma}
\label{th:inverse1}
Every weak group object in a weak 2-category, in particular a stacky Lie group, satisfies the following properties:
\begin{itemize}
 \item[(i)] The inverse is unique up to 2-isomorphism.
 \item[(ii)] A left weak inverse is always a right weak inverse.
 \item[(iii)] The inverse is involutive up to 2-isomorphism.
 \item[(iv)] The weak inverse is a weak antihomomorphism of the group.
\end{itemize}

\end{Lemma}
\begin{proof}
Because the statements are all up to 2-isomorphism, we can consider the associated strict group object in the category of groupoids and isomorphism classes of bibundles. The proof is completely analogous to the case of ordinary set-theoretic groups but still instructive in terms of the graphic notation: First we show that the inverse is unique up to 2-isomorphism. Let $s$ be another inverse. Then:
\begin{equation*}
\begin{xyc}
  \Iden(0,1)\Morp(0,-1){s}\Iden(0,-2)
\end{xyc}
\quad = \quad
\begin{xyc}
\Prod(0,5)
\Iden(-1,3)\Iden(1,3)
\Prod(-1,2) \Iden(1,2)
\Anti(-2,0) \Iden(0,0)\Morp(1,0){s}
\Copr(-1,-1) \Iden(1,-1) \Iden(1,-2)
\Copr(0,-3)
\end{xyc}
\quad = \quad
\begin{xyc}
\Prod(-1,5)
\Prod(0,3) \Iden(-2,3)\Iden(-2,2)
\Iden(-2,1) \Iden(-1,1)
\Anti(-2,0) \Iden(-1,0)
\Copr(0,-1) \Morp(1,0){s} 
\Iden(-2,-1)\Iden(-2,-2)
\Copr(-1,-3)
\end{xyc}
\quad = \quad
\begin{xyc}
  \Iden(0,0)\Anti(0,-1)\Iden(0,-2)
\end{xyc} \quad.
\end{equation*}
A left inverse is involutive up to 2-isomorphism:
\begin{equation*}
\begin{xyc}
  \Iden(0,1)\Iden(0,0)\Iden(0,-1)\Iden(0,-2)
\end{xyc}
\quad = \quad
\begin{xyc}
\Prod(0,5)
\Anti(-1,3) \Iden(1,3)
\Prod(-1,2) \Iden(1,1)\Iden(1,2)
\Anti(-2,0) \Iden(0,0)\Iden(1,0)
\Copr(-1,-1)\Iden(1,-1)\Iden(1,-2)
\Copr(0,-3)
\end{xyc}
\quad = \quad
\begin{xyc}
\Prod(0,5)
\Prod(-1,3) \Iden(1,3)\Iden(1,2)
\Anti(-2,1) \Anti(0,1) \Iden(1,1)
\Anti(-2,0) \Iden(0,0) \Iden(1,0)
\Copr(-1,-1)\Iden(1,-1)\Iden(1,-2)
\Copr(0,-3)
\end{xyc}
\quad = \quad
\begin{xyc}
\Prod(-1,5)
\Prod(0,3) \Iden(-2,3)\Iden(-2,2)
\Anti(-2,1) \Anti(-1,1) \Iden(1,1)
\Anti(-2,0) \Iden(-1,0) \Iden(1,0)
\Copr(0,-1)\Iden(-2,-1)\Iden(-2,-2)
\Copr(-1,-3)
\end{xyc}
\quad = \quad
\begin{xyc}
  \Iden(0,1)\Anti(0,0)\Anti(0,-1)\Iden(0,-2)
\end{xyc} \quad.
\end{equation*}
Using this we can show that a left inverse is also right inverse:
\begin{equation*}
\begin{xyc}
  \Prod(0,2)
  \Iden(-1,0)\Anti(1,0)
  \Copr(0,-1)
\end{xyc}
\quad = \quad
\begin{xyc}
  \Prod(0,3)
  \Anti(-1,1)\Anti(1,1)
  \Anti(-1,0)\Iden(1,0)
  \Copr(0,-1)
\end{xyc}
\quad = \quad
\begin{xyc}
  \Anti(0,3)
  \Prod(0,2)
  \Anti(-1,0)\Iden(1,0)
  \Copr(0,-1)
\end{xyc}
\quad = \quad
\begin{xyc}
  \Anti(0,3)
  \Unit(0,2)
  \Coun(0,-1)
\end{xyc}
\quad = \quad
\begin{xyc}
  \Unit(0,2)
  \Coun(0,0)
\end{xyc}
\quad.
\end{equation*}
The proof that the inverse is an antihomomorphism of the stacky group is analogous but lengthy, so we omit it here.
\end{proof}

In group theoretical language property (iv) states that the inverse is an antihomomorphism of the stacky group. The conclusion is that a weak bimonoid object $G$ together with an inverse is also a weak Hopf object in the category $\LieGrpdPrBibu$. This terminology is in analogy to a Hopf algebra which is a Hopf object in the category of vector spaces. The concept of Hopf objects is very natural and has been studied before \cite{Kuperberg:Involutory}.

So far we have proved properties of the inverse which are analogous to properties for ordinary groups. The following characterization of the stacky group inverse does not have a direct counterpart for set-theoretic groups:

\begin{Theorem}
\label{th:HopfishInverse}
Let $(G, \mu, e)$ be a stacky Lie monoid. Consider the bibundle defined by
\begin{equation*}
  s :=
  \bigl( (\Cv \times \Id \times e) \circ 
  (\Id \times \mu \times \Id)\bigr)  \circ
  (\Id \times \Ev) \,,
\end{equation*}
which we call the \emph{preinverse}. The following are equivalent:
\begin{itemize}
 \item[(i)] The preinverse is a weak isomorphism.
 \item[(ii)] $(G,\mu,e)$ is a stacky Lie group.
\end{itemize}
If these conditions are satisfied then the preinverse is an inverse.
\end{Theorem}

\begin{proof}
 The isomorphism class of the preinverse is represented graphically by
\begin{equation*}
  [s] \quad = \quad
\begin{xyc}
  \Morp(0,-0.5){s}
\end{xyc}
\quad = \quad
  \begin{xyc}
  \Eval(-2,2)
  \Iden(-3,0) \Iden(1,1)
  \Copr(0,0)\Coev(1,-2) \Coun(2,-1)
  \qquad\qquad.
\end{xyc}
\end{equation*}
(i) $\Rightarrow$ (ii): Assume $s$ is a weak isomorphism. Then we have the following isomorphisms of bibundles:
\begin{equation*}
\begin{split}
\begin{xyc}
  \Copr(0,-1) \Morp(-1,0){s} \Iden(1,1)\Iden(1,0) \Prod(0,3)
\end{xyc}
\quad &= \quad
\begin{xyc}
  \Copr(0,-1) \Morp(-1,0){s} \Iden(1,1)\Iden(1,0) 
  \Flip(0,3)\Prod(0,5)
\end{xyc}
\quad = \quad
\begin{xyc}
 \Copr(-2,0)\Coev(2,0)
 \Iden(-3,1)\Iden(-3,2) \Flip(0,2) \Morp(3,1){s}
 \Eval(-2,4)\Prod(2,4)
\end{xyc}
\quad = \quad
\begin{xyc}
 \Copr(-2,0)\Copr(2,0)\Coev(3,-2)\Coun(4,-1)
 \Iden(-3,1)\Iden(-3,2) \Flip(0,2) \Iden(3,1)\Iden(3,2)
 \Prod(-2,4)\Prod(2,4)\Eval(-3,6)\Unit(-4,4)
\end{xyc} \\ 
&= \quad
\begin{xyc}
 \Eval(-2,4)\Unit(-3,2) 
 \Prod(0,0)\Copr(0,2)
 \Coev(2,-2)\Coun(3,-1)
\end{xyc}
\quad = \quad
\begin{xyc}
 \Eval(-2,4)\Unit(-3,2) 
 \Copr(0,2)\Coev(1,0)\Coun(2,1)
 \Coun(3,1)
\end{xyc}
\quad = \quad
\begin{xyc}
 \Morp(0,2){s}
 \Unit(0,1) 
 \Coun(0,-1)
\end{xyc}
\quad = \quad
\begin{xyc}
 \Unit(0,2) \Coun(0,0)
\end{xyc}
\quad,
\end{split}
\end{equation*}
where we have used in the last step that by Proposition~\ref{th:isoprincipal1} the bibundle $s$ is right transitive, and in the second last step we have used that,
\begin{equation*}
\begin{split}
\begin{xyc}
 \Prod(0,0)\Coev(2,-2)\Coun(3,-1)
\end{xyc}
\quad &= \quad
\begin{xyc}
 \Eval(-1,2)\Iden(-2,0)\Iden(-2,-1)
 \Coev(-3,-2)
 \Prod(0,0)
 \Coev(2,-2)\Coun(3,-1)
\end{xyc}
\quad = \quad
\begin{xyc}
 \Eval(0,2)\Iden(1,0)\Iden(1,-1)
 \Coev(-3,-2)
 \Prod(-1,0)
 \Coev(2,-2)\Coun(3,-1)
\end{xyc} 
\quad = \quad
\begin{xyc}
 \Coev(-3,-2)\Coun(-1,1)\Prod(-1,0)
\end{xyc}
\quad = \quad
\begin{xyc}
 \Coev(-3,-2)\Coun(-2,-1)\Coun(0,-1)
\end{xyc}\quad.
\end{split}
\end{equation*}
So far we have shown that $s$ is a left inverse. That $s$ is also a right inverse follows from Lemma~\ref{th:inverse1}.

(ii) $\Rightarrow$ (i): Assume that $G$ is a group object. Consider the isomorphism classes of bibundles
\begin{equation*}
[M] \quad := \quad
\begin{xyc}
  \Iden(-2,1)\Iden(-2,0)\Prod(1,1)
  \Copr(-1,-1) \Iden(2,-1)\Iden(2,-2) 
\end{xyc}
\qquad\text{and}\qquad
[N] \quad := \quad
\begin{xyc}
  \Iden(-2,1)\Iden(-2,0)\Prod(1,1)
  \Iden(-2,-1)\Anti(0,-1)\Iden(2,-1)
  \Copr(-1,-2) \Iden(2,-2)\Iden(2,-3) 
\end{xyc}
\quad.
\end{equation*}
The composition $[M] \circ [N]$ is isomorphic to
\begin{equation*}
\begin{xyc}
  \Iden(-2,1)\Iden(-2,0)\Prod(1,1)
  \Copr(-1,-1) \Iden(2,-1)
  \Prod(2,-1)
  \Iden(-1,-3)\Anti(1,-3)\Iden(3,-3)
  \Copr(0,-4) \Iden(3,-4)\Iden(3,-5) 
\end{xyc}
\quad=\quad
\begin{xyc}
  \Prod(1,4)
  \Prod(0,2)
  \Iden(-2,4)\Iden(-2,3)\Iden(-2,2)\Iden(-2,1)\Iden(-2,0)\Iden(-2,-1)\Iden(-2,-2)
  \Iden(-1,0)\Anti(1,0)
  \Iden(2,1)\Iden(2,2)\Iden(2,0)\Iden(2,-1)\Iden(2,-2)\Iden(2,-3)\Iden(2,-4)
  \Copr(0,-1)
  \Copr(-1,-3)
\end{xyc}
\quad=\quad
\begin{xyc}
  \Prod(1,4)
  \Unit(0,2)
  \Iden(-2,4)\Iden(-2,3)\Iden(-2,2)\Iden(-2,1)\Iden(-2,0)\Iden(-2,-1)\Iden(-2,-2)
  \Iden(2,1)\Iden(2,2)\Iden(2,0)\Iden(2,-1)\Iden(2,-2)\Iden(2,-3)\Iden(2,-4)
  \Coun(0,-2)
  \Copr(-1,-3)
\end{xyc}
\quad=\quad
\begin{xyc}
  \Iden(0,0)\Iden(0,1)\Iden(0,2)
  \Iden(2,0)\Iden(2,1)\Iden(2,2)
\end{xyc}
\quad,
\end{equation*}
that is, $[M] \circ [N] = [\Id_{G\times G}]$. Similarly, it is shown that $[N] \circ [M] = [\Id_{G\times G}]$. So $[N]$ and $[M]$ are mutually inverse. Corollary~\ref{th:isoprincipal2} now implies that $[N] = [M^{\op}]$,
\begin{equation*}
\begin{xyc}
  \Iden(-2,1)\Iden(-2,0)\Prod(1,1)
  \Iden(-2,-1)\Anti(0,-1)\Iden(2,-1)
  \Copr(-1,-2) \Iden(2,-2)\Iden(2,-3) 
\end{xyc}
\quad = \quad
\begin{xyc}
  \Eval(-3,3) \Eval(2,3)
  \Iden(-2,1)\Iden(-2,0)\Prod(1,1)
  \Copr(-1,-1) \Iden(2,-1)\Iden(2,-2) 
  \Coev(-2,-3) \Coev(3,-3)
\end{xyc}
\quad.
\end{equation*}
Composing this relation with $[e] \times [\Id]$ from the left and $[\Id] \times [\varepsilon]$ from the right we obtain
\begin{equation*}
\begin{xyc}
  \Iden(0,1)
  \Anti(0,0)
  \Iden(0,-1)
\end{xyc}
\quad = \quad
\begin{xyc}
  \Coun(-2,1)\Iden(-2,0)\Prod(1,1)
  \Iden(-2,-1)\Anti(0,-1)\Iden(2,-1)
  \Copr(-1,-2) \Iden(2,-2)\Unit(2,-3)
\end{xyc}
\quad = \quad
\begin{xyc}
  \Eval(-3,3) \Eval(2,3)
  \Unit(3,1)
  \Iden(-2,1)\Iden(-2,0)\Prod(1,1)
  \Copr(-1,-1) \Iden(2,-1)\Iden(2,-2) 
  \Coev(-2,-3) \Coev(3,-3)
  \Coun(-3,-2)
\end{xyc}
\quad = \quad
\begin{xyc}
  \Eval(-2,2)
  \Iden(-3,0) \Iden(1,1)
  \Copr(0,0)\Coev(-1,-2) \Coun(-2,-1)
\end{xyc}
\quad = \quad
\begin{xyc}
  \Eval(-2,2)
  \Iden(-3,0) \Iden(1,1)
  \Copr(0,0)\Coev(1,-2) \Coun(2,-1)
\end{xyc}
\quad,
\end{equation*}
which shows that $s$ is isomorphic to the inverse $i$ and hence a weak isomorphism.
\end{proof}

The analog of this statement for set theoretic groups is the following: Let $\mu: G \times G \rightarrow G$ and $e: 1 \rightarrow G$ be the morphisms of a semi-group structure on the set $G$. The pull-back $(G \times G) \times_{G}^{\mu,e} 1$ can be identified with the set of mutually inverse elements $S := \{(g,h) \in G \times G \,|\, gh = e \}$. The semi-group $G$ is a group if and only if $S$ is the graph of a map. If it is, then the map is automatically a bijection. The remarkable difference in the stacky case is the existence of the preinverse, which is a morphism in the category of groupoids and bibundles. The stacky preinverse is the geometric analog of the algebraic concept of the hopfish preantipode \cite{BlohmannTangWeinstein:Hopfish} \cite{TangWeinsteinZhu:Hopfish} and has already appeared in \cite{Zhu:n-groupoids}.

\subsection{The stacky Lie group of a reduced torus foliation}

A Lie 2-group is a group object in the category of Lie groupoids and differentiable homomorphisms of groupoids. Bundlization of the structure morphisms of a Lie 2-group yields the structure bibundles of a stacky Lie group. Therefore, Lie 2-groups are a special case of stacky Lie groups and a natural source of examples. While not every stacky Lie group is a strict Lie 2-group, it can be shown that every stacky Lie group is a weak Lie 2-group up to Morita equivalence \cite{Zhu:LieII} \cite{Zhu:n-groupoids}.

There is a 1-to-1 correspondence between strict 2-groups and crossed modules \cite{BaezLauda:2-Groups} \cite{Breen:Bitorseurs} \cite{Noohi:2-groupoids}. An obvious example of a crossed module is given by the inclusion $\phi: A \hookrightarrow B$ of a normal subgroup with the adjoint action $b \cdot a = b a b^{-1}$. For example if $A=B$ is a group we get a group structure on the groupoid of the left regular action of $A$ on itself. The coarse moduli space of the corresponding stack is just a point. 

For a more interesting example we return to our guiding example. Consider the plane $\bbC$, with its standard $S^1 = U(1)$ symmetry. The action of $\bbZ$ on $\bbC$ by rotation by a finite angle $k \cdot z = e^{i\lambda k} z$ commutes with the $S^1$-action. In fact, if $\lambda/2\pi$ is irrational we can faithfully embed $\bbZ$ as dense normal subgroup $e^{i \lambda \bbZ} \hookrightarrow S^1$. This yields a crossed module and, hence, a strict 2-group. Let us describe the 2-group structure explicitly \cite{BlohmannTangWeinstein:Hopfish}:

As a set, the 2-group is given by $\gT_1 = A \times B = \bbZ \times S^1$ which is a groupoid over $\gT_0 = B = S^1$. The left and right moment maps are given by $l_{\gT}(k,\theta) = \theta + \lambda k$ and $r_{\gT}(k,\theta) = \theta$, where the argument $\theta$ is shorthand for $e^{i \theta}$. If we think of a groupoid element as an arrow pointing from an element $\theta$ of $S^1$ to its image under the action of $k \in \bbZ$, then $r$ is the source map and $l$ the target. When $r(k_1,\theta_1) = l(k_2,\theta_2)$, that is, when $\theta_1 = \theta_2 + \lambda k_2$, the groupoid product is given by $(k_1,\theta_1)(k_2,\theta_2) = ( k_1 + k_2, \theta_2 )$. The identity bisection is given by the embedding $S^1 \hookrightarrow S^1 \times \bbZ$, $\theta \mapsto (0,\theta)$.  The groupoid inverse is given by $(k,\theta)\mapsto (-k, \theta -\lambda k)$. With the product differentiable structure $\gT$ is a smooth, \'etale groupoid. Note, however, that the anchor map $(l,r): \gT_1 \rightarrow \gT_0 \times \gT_0$ is not proper, so $\gT_1$ is not the presentation of an orbifold.

The group structure on $\gT = \bbZ \times S^1$ is given by the product group structure, $\mu\bigl((k_1,\theta_1),(k_2,\theta_2)\bigr) = (k_1 + k_2,\theta_1 + \theta_2)$, $e = (0,0)$, $i(k,\theta) = (-k,-\theta)$. By bundlization of the multiplication we get the bibundle
\begin{equation*}
\begin{split}
  \Bibu{\mu} 
  &:= (\gT_0 \times \gT_0) \times_{\gT_0}^{\mu,l} \gT_1
  = (S^1 \times S^1) \times_{S_1} (\bbZ \times S^1) \cong S^1 \times S^1 \times \bbZ \,,
\end{split}
\end{equation*}
with moment maps $l_{\Bibu{\mu}}(\theta_1,\theta_2, k) = (\theta_1,\theta_2)$ and $r_{\Bibu{\mu}}(\theta_1,\theta_2, k) = \theta_1 + \theta_2 + \lambda k$ for $(\theta_1,\theta_2,k) \in S^1 \times S^1 \times \bbZ$. The left action of $(l_1, \phi_1, l_2, \phi_2) \in \gT \times \gT$ and right action of $(l,\phi) \in \gT$ is
\begin{equation*}
\begin{aligned}
 (l_1,\phi_1,l_2,\phi_2) \cdot (\theta_1, \theta_2, k) 
  &= (\theta_1 + \lambda l_1, \theta_2 + \lambda l_2, k - l_1 - l_2) \,,\\
 (\theta_1, \theta_2, k) \cdot (l, \phi)
  &= (\theta_1, \theta_2, k - l) \,,
\end{aligned}
\end{equation*}
which is defined if $\phi_1 = \theta_1$, $\phi_2 = \theta_2$, and $\theta_1 + \theta_2 + \lambda k = \phi + \lambda l$.

The bundlization of the group identity is
\begin{equation*}
  \Bibu{e} := 1 \times_{\gT_0}^{e,l} \gT_1
  = 1 \times_{S_1} (\bbZ \times S^1) \cong \bbZ \,,
\end{equation*}
with moment maps $l_{\Bibu{e}}(k) = 1$ and $r_{\Bibu{e}}(k) = - \lambda k$ for $k \in \bbZ$. The left groupoid action is trivial, whereas the right $\gT$-action is given for $(l,\phi) \in \gT_1$ by
\begin{equation*}
  k \cdot (l,\phi) = k + l \,,
\end{equation*}
which is defined if $-\lambda k = \phi + \lambda l$. Finally, the bundlization of the inverse is given by
\begin{equation*}
  \Bibu{i} := \gT_0 \times_{\gT_0}^{i,l} \gT_1
  = S_1 \times_{S_1} (\bbZ \times S^1) \cong S^1 \times \bbZ \,,
\end{equation*}
with moment maps $l_{\Bibu{i}}(\theta,k) = \theta$ and $r_{\Bibu{i}}(\theta,k) = -\theta - \lambda k$ for $(\theta,k) \in S^1 \times \bbZ$. The left and right action of $(l,\phi) \in \gT_1$ is
\begin{equation*}
  (l,\phi) \cdot (\theta,k) = (\theta + \lambda l, k - l) \,, \qquad
  (\theta, k) \cdot (l,\phi) = (\theta, k + l) \,,
\end{equation*}
which is defined if $\phi =  \theta$ and if $-\theta - \lambda k = \phi + \lambda l$, respectively.

With the product differentiable structure these three bibundles are smooth right principal bibundles. By functoriality of bundlizations, they satisfy the axioms for the structure bibundles of a stacky Lie group. The interesting fact to note here is that for irrational $\lambda/2\pi$ the coarse moduli space of the groupoid $\gT$ inherits the trivial quotient topology, whereas for rational values it is diffeomorphic to $S^1$. But if we look at the family of groupoids $\gT$ parametrized by $\lambda$ we have a smooth deformation of stacky Lie groups.

\subsection*{Acknowledgements}

I am indebted to Alan Weinstein for uncountable discussions and a large number of improvements he has brought to this work. I owe to him the idea of linking categories and linking groupoids. I would like to thank the following colleagues for helpful discussions of aspects of this work: Larry Breen, Kai Behrend, Andr\'e Haefliger, Andrew Kresch, Ieke Moerdijk, Martin Olsson, Chris Schommer-Pries, Marc Rieffel, Xiang Tang, Peter Teichner, and Chenchang Zhu. I would also like to thank Ping Xu and the Centre Emile Borel at the Institut Henri Poincar\'e for their kind hospitality. This work was supported by a Marie Curie Fellowship of the European Union

\providecommand{\bysame}{\leavevmode\hbox to3em{\hrulefill}\thinspace}
\providecommand{\href}[2]{#2}


\begin{thebibliography}{10}

\bibitem{BaezLauda:2-Groups}
John~C. Baez and Aaron~D. Lauda, \emph{Higher-dimensional algebra. {V}.
  2-groups}, Theory Appl. Categ. \textbf{12} (2004), 423--491 (electronic).

\bibitem{BehrendXu:Differentiable}
Kai Behrend and Ping Xu, \emph{{Differentiable Stacks and Gerbes}},
  \mbox{arXiv:math.DG/0605694}.

\bibitem{BehrendXu:S1-bundles}
Kai Behrend and Ping Xu, \emph{{$S\sp 1$}-bundles and gerbes over
  differentiable stacks}, C. R. Math. Acad. Sci. Paris \textbf{336} (2003),
  no.~2, 163--168.

\bibitem{BlohmannTangWeinstein:Hopfish}
Christian Blohmann, Xiang Tang, and Alan Weinstein, \emph{{Hopfish structure
  and modules over irrational rotation algebras}}
  \mbox{arXiv:math.QA/0604405}, to appear in Contemp.~Math.

\bibitem{BlohmannWeinstein:Group-like}
Christian Blohmann and Alan Weinstein, \emph{{Group-like objects in Poisson
  geometry and algebra}}, Contemp.~Math.~\textbf{450} (2008), 15--39, \mbox{arXiv:math.SG/0701499}.

\bibitem{Blohmann:Higher}
  Christian Blohmann, in preparation.

\bibitem{Borceux:Handbook1}
Francis Borceux, \emph{Handbook of categorical algebra. 1. Basic category theory.}, Encyclopedia of Mathematics and its Applications \textbf{50}, Cambridge University Press, Cambridge, 1994.

\bibitem{Breen:Bitorseurs}
Lawrence Breen, \emph{Bitorseurs et cohomologie non ab\'elienne}, The
  Grothendieck Festschrift, Vol.\ I, Progr. Math., vol.~86, Birkh\"auser
  Boston, Boston, MA, 1990, pp.~401--476.

\bibitem{BrownGreenRieffel:Stable}
Lawrence~G. Brown, Philip Green, and Marc~A. Rieffel, \emph{Stable isomorphism
  and strong {M}orita equivalence of {$C\sp*$}-algebras}, Pacific J. Math.
  \textbf{71} (1977), no.~2, 349--363.

\bibitem{CrainicFernandes:Integrability}
Marius Crainic and Rui~Loja Fernandes, \emph{Integrability of {L}ie brackets},
  Ann. of Math. (2) \textbf{157} (2003), no.~2, 575--620.

\bibitem{DeligneMumford:irreducibility}
P.~Deligne and D.~Mumford, \emph{The irreducibility of the space of curves of
  given genus}, Inst. Hautes \'Etudes Sci. Publ. Math. (1969), no.~36, 75--109.

\bibitem{Freed:Higher}
Daniel~S. Freed, \emph{Higher algebraic structures and quantization}, Comm.
  Math. Phys. \textbf{159} (1994), no.~2, 343--398.

\bibitem{Giraud:Cohomologie}
Jean Giraud, \emph{Cohomologie non ab\'elienne}, Springer-Verlag, Berlin, 1971,
  Die Grundlehren der mathematischen Wissenschaften, Band 179.

\bibitem{Heinloth:Notes}
J.~Heinloth, \emph{Notes on differentiable stacks}, Mathematisches Institut,
  Georg-August-Universit\"at G\"ottingen: Seminars Winter Term 2004/2005,
  Universit\"atsdrucke G\"ottingen, G\"ottingen, 2005, pp.~1--32.

\bibitem{HilsumSkandalis:Morphismes}
Michel Hilsum and Georges Skandalis, \emph{Morphismes {$K$}-orient\'es
  d'espaces de feuilles et fonctorialit\'e en th\'eorie de {K}asparov
  (d'apr\`es une conjecture d'{A}. {C}onnes)}, Ann. Sci. \'Ecole Norm. Sup. (4)
  \textbf{20} (1987), no.~3, 325--390.

\bibitem{Joyal:Quasi}
Andr{\'e} Joyal, \emph{Quasi-categories and {K}an complexes}, J.~Pure~Appl.~Alg. \textbf{175} (2002), 207--222.

\bibitem{Kuperberg:Involutory}
Greg Kuperberg, \emph{Involutory {H}opf algebras and {$3$}-manifold
  invariants}, Internat. J. Math. \textbf{2} (1991), no.~1, 41--66.

\bibitem{Landsman:Quantized}
N.~P. Landsman, \emph{Quantized reduction as a tensor product}, Quantization of
  singular symplectic quotients, Progr. Math., vol. 198, Birkh\"auser, Basel,
  2001, pp.~137--180.

\bibitem{Landsman:Lie_groupoids}
\bysame, \emph{Lie groupoids and {L}ie algebroids in physics and noncommutative
  geometry}, J. Geom. Phys. \textbf{56} (2006), no.~1, 24--54.

\bibitem{LupercioUribe:Gerbes}
Ernesto Lupercio and Bernardo Uribe, \emph{Gerbes over orbifolds and twisted
  {$K$}-theory}, Comm. Math. Phys. \textbf{245} (2004), no.~3, 449--489.

\bibitem{Lurie:Topos}
Jacob Lurie, \emph{Higher topos theory}, \mbox{arXiv:math.CT/0608040}

\bibitem{Metzler:Topological}
David Metzler, \emph{{Topological and Smooth Stacks}},
  \mbox{arXiv:math.DG/0306176}.

\bibitem{MoerdijkPronk:Orbifolds}
I.~Moerdijk and D.~A. Pronk, \emph{Orbifolds, sheaves and groupoids},
  $K$-Theory \textbf{12} (1997), no.~1, 3--21.

\bibitem{Moerdijk:Orbifolds}
Ieke Moerdijk, \emph{Orbifolds as groupoids: an introduction}, Orbifolds in
  mathematics and physics (Madison, WI, 2001), Contemp. Math., vol. 310, Amer.
  Math. Soc., Providence, RI, 2002, pp.~205--222.

\bibitem{Mrcun:Stability}
Janez Mrcun, \emph{{Stability and invariants of Hilsum-Skandalis maps}}, Ph.D.
  thesis, \mbox{arXiv:math.DG/0506484}.

\bibitem{Mrcun:Functoriality}
Janez Mr{\v{c}}un, \emph{Functoriality of the bimodule associated to a
  {H}ilsum-{S}kandalis map}, $K$-Theory \textbf{18} (1999), no.~3, 235--253.

\bibitem{Noohi:2-groupoids}
Behrang Noohi, \emph{{Notes on 2-groupoids, 2-groups and crossed-modules}},
  \mbox{arXiv:math.CT/0512106}.

\bibitem{Pradines:Theorie}
Jean Pradines, \emph{Th\'eorie de {L}ie pour les groupo\"\i des
  diff\'erentiables. {R}elations entre propri\'et\'es locales et globales}, C.
  R. Acad. Sci. Paris S\'er. A-B \textbf{263} (1966), A907--A910.

\bibitem{TangWeinsteinZhu:Hopfish}
Xiang Tang, Alan Weinstein, and Chenchang Zhu, \emph{{Hopfish algebras}},
  \mbox{arXiv:math.QA/0510421}.

\bibitem{ZhuTsegn:Integrating}
Hsian-Hua Tseng and Chenchang Zhu, \emph{Integrating {L}ie algebroids via
  stacks}, Compos. Math. \textbf{142} (2006), no.~1, 251--270.

\bibitem{TuXuLaurent-Gengoux:Twisted}
Jean-Louis Tu, Ping Xu, and Camille Laurent-Gengoux, \emph{Twisted {$K$}-theory
  of differentiable stacks}, Ann. Sci. \'Ecole Norm. Sup. (4) \textbf{37}
  (2004), no.~6, 841--910.

\bibitem{Zhu:LieII}
Chenchang Zhu, \emph{{Lie II theorem for Lie algebroids via stacky Lie
  groupoids}}, \mbox{arXiv:math.DG/0701024}.

\bibitem{Zhu:n-groupoids}
\bysame, \emph{{Lie n-groupoids and stacky Lie groupoids}},
  \mbox{arXiv:math.DG/0609420}.

\end{thebibliography}
\end{document}